\documentclass[12pt,reqno]{amsart}
\usepackage{graphicx}
\usepackage{subfigure}
\usepackage{times}
\usepackage{amsmath,amsfonts, amstext,amssymb,amsbsy,amsopn,amsthm}
\usepackage{dsfont}
\usepackage{esint}
\usepackage{graphicx}   
\usepackage{hyperref}
\usepackage[all]{xy}
\usepackage{color}
\usepackage{tikz}
 \usetikzlibrary{intersections}
\usetikzlibrary{external}
\usepackage{lineno,hyperref}
\allowdisplaybreaks[4]

\newcommand{\la}{\lambda}
\newcommand{\lap}{\mbox{$\bigtriangleup$}}
\newcommand{\flap}{\mbox{$(-\triangle)^s$}}
\newcommand{\norm}[1]{\lVert#1\rVert}
\newcommand{\grad}{\mbox{$\bigtriangledown$}}
\newcommand{\ra}{{\mbox{$\rightarrow$}}}
\newcommand{\be}{\begin{equation}}
\newcommand{\ee}{\end{equation}}
\newcommand{\R}{\mathbb{R}}

\newcommand{\hspn}{ \Sigma_{\lambda}^{-} }
\newcommand{\hsp}{ \Sigma_{\lambda} }
\newcommand{\dR}{\mathds{R}}

\newcommand{\bee}{\begin{equation*}}
\newcommand{\eee}{\end{equation*}}

\newtheorem{mdef}{Definition}

\newtheorem{thm}{Theorem}[section]
\newtheorem{pro}{Proposition}[section]
\newtheorem{lem}{Lemma}[section]
\newtheorem{cor}{Corollary}[section]
\newtheorem{rem}{Remark}[section]

\begin{document}

\title[Methods on  fractional equations]{Methods in studying qualitative properties of fractional equations}

\author{Wenxiong Chen \quad  Yahong Guo \quad  Congming Li* }
\footnotetext[1]{* Corresponding author.}

\begin{abstract}

In this paper, we systematically review a series of effective methods for studying the qualitative properties of solutions to fractional equations. Beginning with the pioneering extension method and the method of moving planes in integral forms, we introduce a variety of direct methods, including the direct method of moving planes, the method of moving spheres, blow-up and rescaling techniques, the sliding method, regularity lifting, and approaches for interior and boundary regularity estimates.

To elucidate the core ideas behind these methods, we employ simple examples that demonstrate how they can be applied to investigate qualitative properties of solutions. We also provide a comparative discussion of their respective strengths and limitations. It is our hope that this paper will serve as a useful handbook for researchers engaged in the study of fractional equations.

\end{abstract}

\maketitle

Mathematics Subject classification (2020): 35R11, 35B50, 35B06, 35A01.
\medskip

Keywords:  The fractional Laplacian; nonlocal elliptic and parabolic equations; direct method of moving planes; method of moving spheres; sliding method; blow-up and rescaling; a priori estimates; qualitative properties, radial symmetry; monotonicity; nonexistence; regularity lifting; interior regularity; boundary regularity.

\bigskip

\pagebreak

\tableofcontents

\pagebreak

\section{Introduction}
\medskip

 The fractional Laplacian is a nonlocal operator defined by the singular integral
\be
(-\lap)^s u (x) \equiv C_{n, s} PV \int_{\mathbb{R}^n}\frac{u(x)-u(y)}{|x-y|^{n+2s}}dy
\label{1.1}
\ee
for any real number $0<s<1$, where PV stands for the Cauchy principal value, and $C_{n,s}$ is a normalization constant.

This is a singular integral involving two singularities: one at the point $x$, the other near $\infty$.

In order the integral to converge in a neighborhood of $x$, we require $u(x)$ to be locally in $C^{1,1}$ (actually, in $C^{s+\epsilon}$ for some $\epsilon >0$ is sufficient).

To ensure the convergence of the singular integral near infinity, we assume that $u$ is a slowly increasing function in the space:
$${\mathcal L}_{2s}=\left\{u \in L_{\rm loc}^1(\mathbb{R}^n) \mid \int_{\mathbb{R}^n}\frac{|u(x)|}{1+|x|^{n+2s}} \, d x <\infty\right\}.$$

 It can be verified that, for each fixed $x$, as $s \to 1,$
$$(-\lap)^su(x) \to - \lap u(x),$$
where $\lap$ is the well-known Laplace operator. That's why we call $(-\lap)^s$ the fractional Laplacian.
\smallskip

Due to its ability to describe nonlocal, anomalous  effects and long range interactions makes the fractional Laplacian a powerful tool across diverse scientific and engineering domains and has numerous applications in Physics and Engineering, Image Processing and Computer Vision, Probability and Stochastic Processes, Mathematical Biology and Ecology, Quantum Mechanics and Relativistic Models, Materials Science and Mechanics, Machine Learning and  Data Science, and Geophysics and Climate Modeling, and so on.
\medskip

Besides the fractional Laplacian, there are other elliptic fractional operators as listed in the following.

\begin{enumerate}

\item  {\em Uniformly elliptic fractional operator}
 $$ (-\Delta)_a^s u(x) = C_{n,s} P. V. \int_{\mathbb{R}^n} \frac{a(x-y)(u(x)-u(y))}{|x-y|^{n+2s}} d y $$
with the function $a(\cdot)$ bounded from above and bounded away from $0$.

\item  {\em Nonlocal Monge-Amp\`{e}re operator}
$$ \mathcal{D}_s u(x) = \inf\left\{ P. V. \int_{\mathbb{R}^n} \frac{u(y)-u(x)}{|A^{-1}(y-x)|^{n+2s}} d y\mid A>0, \det A=1\right\},$$
where $A>0$ means that the square matrix $A$ is positive definite.

To ensure such operators obey maximum principles, one usually consider their sub-family
$$ D_s^\theta = \{ D_s \mid \lambda_{\min} (A) \geq \theta > 0\}$$
with $\lambda_{\min}(A)$ being the minimum eigenvalue of the matrix $A$. Such $D_s^\theta$ is uniformly elliptic.

\item {\em Fully nonlinear nonlocal operators and the fractional p-Laplacian}

$$
F_{\alpha}(u(x)) = C_{n,\alpha} = C_{n,\alpha} PV \int_{\mathbb{R}^n} \frac{G(u(x)-u(y))}{|x-y|^{n+\alpha}} d y.
$$
Here $0< \alpha <2$ and $G$ is a nonlinear function that is at least local Lipschitz continuous with $G(0)=0$.

In particular,

  (i) When $G(\cdot)$ is an identity map, $F_{\alpha}$ becomes the fractional Laplacian $(-\lap)^{\alpha/2}$.

  ii) When $G(t)= |t|^{p-2}t$ and $\alpha = ps $, $F_\alpha$ become the fractional $p$-Laplacian defined by $$(-\lap)^s_p u(x) = C_{n,s} PV \int_{\mathbb{R}^n} \frac{|u(x)-u(y)|^{p-2}(u(x)-u(y))}{|x-y|^{n+ps}}  dy .$$

   \medskip

If $p \geq 2$ is fixed, then as $s \ra 1$, one can show that, for each $x$,
$$(-\lap)^s_p u(x) \ra - div(|\grad u|^{p-2} \grad u(x)) \equiv -\lap_p u(x).$$

\end{enumerate}

\medskip

When studying equations involving fractional operators, one encounters substantial challenges due to their nonlocal nature, which renders many classical analytical methods—originally developed for local operators—inapplicable.

To overcome the difficulties arising from the nonlocality of the fractional Laplacian, Caffarelli and Silvestre \cite{CS} introduced the celebrated extension method, which transforms nonlocal problems involving the fractional Laplacian into equivalent local problems in higher-dimensional spaces. This reformulation allows one to apply the well-developed machinery of classical elliptic theory to analyze and solve the extended equations.

The extension method has been widely and successfully employed to study equations involving the fractional Laplacian, yielding a wealth of deep and influential results (see, for example, \cite{BCPS}, \cite{CZ,CFS,S}, and the references therein).

Another powerful approach to fractional equations is the integral equation method, first developed in \cite{CLO}, \cite{MCL}, and \cite{CL1}. This framework enables the use of tools from integral equation theory—such as the method of moving planes in integral form and regularity lifting techniques—to investigate qualitative properties of solutions.

Both of these methods were designed mainly for equations involving the fractional Laplacians. However there are other nonlocal elliptic operators, such as the above mentioned uniformly elliptic fractional operators, nonlocal Monge-Amp\`{e}re operators, and fully nonlinear nonlocal operators including the fractional p-Laplacians, so far as we know, neither {\em extension methods} nor {\em integral equations methods} have been tailored for equations involving these operators. This gap highlights the need to develop methods that can directly address such nonlocal problems.

In response to this challenge, the authors in \cite{CLL} introduced a direct method of moving planes for fractional equations, which has proven highly effective in deriving qualitative properties of solutions to nonlinear problems involving the fractional Laplacian on various domains. This approach was subsequently modified and extended in \cite{CLg} and \cite{CL2} to handle nonlinear equations governed by fully nonlinear nonlocal operators, including the fractional $p$-Laplacians, leading to a series of interesting results on symmetry, monotonicity, and nonexistence of solutions.

In a related development, the authors in \cite{CLZ} proposed a direct method of moving spheres for fractional equations, which has been successfully applied to explore qualitative properties of solutions in various settings.

Furthermore, in \cite{CLL1}, a direct blow-up and rescaling method was introduced for fractional equations to derive a priori estimates for positive solutions on bounded domains, thereby establishing existence results via topological degree theory.

The sliding method, originally introduced by Berestycki, Caffarelli, and Nirenberg \cite{BCN2, BHM, BN3}, serves as another fundamental tool to obtain monotonicity and one-dimensional symmetry for local elliptic problems. This method has been substantially reformulated in \cite{WuC, WC2} to apply directly to nonlocal equations involving the fractional Laplacian and the fractional $p$-Laplacian. More recently, after essential modifications, it has been successfully employed to establish fractional parabolic analogues of Gibbons’ conjecture \cite{CM1, CW1}.

These newly developed techniques have become powerful and versatile analytical tools for studying qualitative properties of solutions to fractional elliptic equations. They have since been widely adopted and further refined by numerous researchers to address a broad spectrum of nonlocal problems (see, e.g., \cite{BPGQ, CH, CL, CLO2, CLZ1, CM, CW1, DQ, LC, LZ, JLX, ZCCY, ZL} and the references therein).

This paper is organized as follows.

In Section 2, we  provide several different equivalent definitions of the fractional Laplacian on the Euclidian space as well as closed Riemannian manifolds. We also  list two applications of  the {\em Caffarelli and Silvestre's extension method}.

Section 3 presents several typical maximum principles, including {\em maximum principle in a punctured ball} and {\em the maximum principle in unbounded domains}.

In Section 4, we recapitulate the {\em method of moving planes in integral forms}. We show how to derive an equivalent integral equation and highlight the essence
of this method in establishing radial symmetry of solutions.

In Section 5, we summarize various direct methods that have been widely applied to study the fractional equations, including the {\em direct method of moving planes}, {\em method of moving spheres}, {\em of blow-up and rescaling technique}, and {\em the sliding method}. Simple examples are provided to illustrate the ideas in applying these methods to derive qualitative properties of fractional equations, such as symmetry, monotonicity, nonexistence, and uniqueness of solutions. We will compare these methods and discuss their relative strength and limitations.

The last section is devoted to the regularity of solutions. Starting from the classical {\em regularity lifting by contracting or shrinking operators}, we present our most recent results on the interior and boundary regularity estimates for nonlocal elliptic and parabolic equations. Instead of using the global bounds, we employ only the local bounds of nonnegative solutions to control the higher interior norms. This make it possible to obtain a priori estimates of solutions to fractional equations on unbounded domains via blow-up and rescaling technique.

\section{Other equivalent  definitions for fractional Laplacian}
\subsection{The fractional Laplacian on Euclidean space $\mathbb{R}^n$} \,

\medskip

Besides the singular integral definition of the fractional Laplacian \eqref{1.1}, there are other equivalent definitions.

\begin{mdef}[The Fourier transform \cite{S}]Let $u$ be a function in the
Schwartz class  $\mathcal{S}(\mathbb{R}^n), 0<s<1,$ the fractional Laplacian $(-\Delta)^s$ is defined via Fourier transform as
\[\widehat{(-\Delta)^s} u(\xi)=|\xi|^{2s}\widehat{u}(\xi).\]
\end{mdef}

\begin{mdef}[The heat kernel \cite{B,ST1}] For $u\in C_{loc}^{2s+\varepsilon}(\mathbb{R}^n)\cap {\mathcal L}_{2s}$,  the fractional Laplacian $(-\Delta)^s$ can also be introduced  by \[(-\Delta)^su(x)
=\frac{1}{\Gamma(-s)}\int_0^{\infty}(e^{t\Delta}u(x)-u(x))\frac{dt}{t^{1+s}},\]
where $(e^{t\Delta})_{t>0}$ denotes the heat semigroup generated by $-\Delta$, that is,
\[e^{t\Delta}u(x)=\frac{1}{(4\pi t)^{\frac{n}{2}}}\int_{\mathbb{R}^n}e^{\frac{-|x-z|^2}{4t}}{u(z)}dz.\]

\end{mdef}

\begin{mdef}[The extention \cite{CS}]
 For a function $u:\mathbb{R}^n \ra \mathbb{R}$, one considers its extension $U:\mathbb{R}^n\times[0, \infty) \ra \mathbb{R}$ that satisfies
\begin{equation}\label{Etr}\left\{\begin{array}{ll}
div(y^{1-2s} \nabla U)=0, & (x,y) \in \mathbb{R}^n\times[0, \infty),\\
U(x, 0) = u(x).
\end{array}
\right.
\end{equation}
Then
\begin{equation}\label{rpt}(-\lap)^{s}u (x) = - C_{n,s} \displaystyle\lim_{y \ra 0^+}
y^{1-2s} \frac{\partial U}{\partial y},  \;\; x \in \mathbb{R}^n.\end{equation}
\end{mdef}
\subsection{The fractional Laplacian on a closed manifold $(M, g)$}\,
\smallskip

Taking inspiration from the case of $\mathbb{R}^n$, Caselli,  Florit-Simon and Serra \cite{CFS} introduce several equivalent definitions for the fractional Laplacian $(-\Delta)^{s}$ on a closed (i.e. compact and without boundary) Riemannian manifold $(M, g)$ of dimension $n$, with $s \in (0, 1)$.

\subsubsection{Spectral and singular integral definitions}

The fractional Laplacian $(-\Delta)^{s}$ can be defined as the $s$-th power (in the sense of spectral theory) of the usual Laplace--Beltrami operator on a Riemannian manifold, through Bochner's subordination. 



\begin{mdef}[Spectral definition\cite{CFS}]
Let $s \in (0, 1)$. The fractional Laplacian $(-\Delta)^{s}$ is the operator that acts on smooth functions $u$ by

\[
(-\Delta)^{s} u = \frac{1}{\Gamma(-s)} \int_{0}^{\infty} (e^{t\Delta} u - u) \frac{dt}{t^{1+s}}.\tag{SP}
\]

Here, the expression $e^{t\Delta} u$ is to be understood as the solution of the heat equation on $M$ at time $t$ and with initial datum $u$. 
\end{mdef}

\begin{rem}
 The above definition is easily shown to be identical to
 \[
(-\Delta)^{s} u = \sum_{k=1}^\infty \lambda_k^{s} \langle u, \phi_k \rangle_{L^2(M)} \phi_k.
\]
Here, $\{\phi_k\}_{k=1}^\infty$ is an $L^2(M)$ orthonormal basis of eigenfunctions for $(-\Delta)$ with eigenvalues
\[
0 < \lambda_1 < \lambda_2 \leq \ldots \leq \lambda_k \xrightarrow{k \to \infty} + \infty
\]
\end{rem}

The second definition for the fractional Laplacian, which is closely related and equivalent to the spectral definition, expresses it as a singular integral.
\begin{mdef}[Singular integral definition \cite{CFS}]
The fractional Laplacian $(-\Delta)^{s}$ of order (of differentiation) $s \in (0, 1)$ is the operator that acts on a regular function $u$ by
\[
(-\Delta)^{s} u(p) = \mathrm{p.v.} \int_M (u(p) - u(q)) K_s(p, q)  dV_q \tag{SI}
\]
\[
:= \lim_{\epsilon \to 0} \int_M (u(p) - u(q)) K_s^\epsilon (p, q)  dV_q.
\]
Here $K_s(p, q) : M \times M \to \mathbb{R}$ denotes the singular kernel given by
\[
K_s(p, q) = \frac{s}{\Gamma(1-s)} \int_0^\infty H_M(p, q, t) \frac{dt}{t^{1+s}}
\]
where $H_M : M \times M \times (0, \infty) \to \mathbb{R}$ denotes the usual heat kernel on $M$, and $K_s^\epsilon (p, q)$ is the natural regularization
\[
K_s^\epsilon (p, q) = \frac{s}{\Gamma(1-s)} \int_0^\infty H_M(p, q, t) e^{- \epsilon^2 / 4t} \frac{dt}{t^{1+s}}.
\]
\end{mdef}

\begin{rem}
If the compact manifold $M$ is replaced by the Euclidean space $\mathbb{R}^n$ then
\begin{align*}
K_s(x, y) &= \frac{s}{\Gamma(1 - s)} \int_0^\infty H_{\mathbb{R}^n}(x, y, t) \frac{dt}{t^{1 + s}} \\
&= \frac{s}{\Gamma(1 - s)} \int_0^\infty \left( \frac{1}{(4\pi t)^{\frac{n}{2}}} e^{-\frac{(x - y)^2}{4t}} \right) \frac{dt}{t^{1 + s}} = \frac{\alpha_{n,s}}{|x - y|^{n + 2s}},
\end{align*}
Hence we recover the usual form of the fractional Laplacian on $\mathbb{R}^n$.
\end{rem}
\subsubsection{Extension definition}

It was proved by Stinga \cite{S} that the unique solution to \eqref{Etr} verifying \eqref{rpt} admits the explicit representation
\begin{equation}\label{rep}
U(p, z) = \frac{z^{2s}}{2^{2s} \Gamma(s)} \int_0^\infty e^{t\Delta} u(p) e^{-\frac{z^2}{4t}} \frac{dt}{t^{1+s}},
\end{equation}
which expresses \( U \) in terms of the solution to the heat equation \( e^{t\Delta} u \) (and thus makes sense also on a manifold).

\begin{mdef}[Extension definition\cite{CFS}]
Let \((M^n, g)\) be a closed Riemannian manifold, let \( s \in (0, 1) \) and \( u : M \to \mathbb{R} \) be smooth. Consider the product manifold \( \tilde{M} = M \times (0, +\infty) \) endowed with the natural product metric.\footnote{That is, the metric defined by \( \tilde{g}((\xi_1, z_1), (\xi_2, z_2)) = g(\xi_1, \xi_2) + z_1 z_2 \), and where \( \tilde{\mathrm{div}} \) and \( \tilde{\nabla} \) denote the divergence and Riemannian gradient with respect to this product metric respectively.} Then, there is a unique solution \( U : M \times (0, \infty) \to \mathbb{R} \) among functions in \( \tilde{H}^1(\tilde{M}) \) to
\[
\begin{cases}
\tilde{\mathrm{div}}(z^{1-2s}\tilde{\nabla}U) = 0 \quad \text{in} \quad \tilde{M}, \\
U(p, 0) = u(p) \quad \text{for} \quad p \in \partial \tilde{M} = M,
\end{cases}
\]
given by \eqref{rep}, and it satisfies
\[
[u]_{H^{s}(M)}^2 = 2\beta_s \int_{\tilde{M}} |\tilde{\nabla}U|^2 z^{1-2s}  dV  dz,
\]
where \([u]_{H^{s}(M)}^2\) is defined by
\[[u]_{H^{s}(M)}^2:=\int\int\limits_{M\times M}(u(p)-u(q))^2K_s(p,q)dV_pdV_q,\] and
\[
\beta_s = \frac{2^{2s-1} \Gamma(s)}{\Gamma(1-s)}.
\]
Moreover,
\[
\lim_{z \to 0^+} z^{1-2s} \frac{\partial U}{\partial z}(p, z) = -\beta_s^{-1}(-\Delta)^{s} u(p),
\]
where the fractional Laplacian on the right-hand side is defined by either (SP) or  ({SI}).
\end{mdef}
Here, the weighted Sobolev space on the manifold
\[
\tilde{H}^1(\tilde{M}) = \tilde{H}^1({M}\times (0, \infty))=H^1({M} \times (0, \infty), z^{1 - 2s}  dV  dz)
\]
is the completion of \( C_c^\infty(M \times [0, \infty)) \) with the norm
\[
\|U\|_{\tilde{H}^1}^2 := \|TU\|_{L^2(M)}^2 + \|\tilde{\nabla}U\|_{L^2(\tilde{M}, z^{1-2s}  dV  dz)},
\]
where \( TU = U(\cdot, 0) \) is the trace of \( U \) and \( \tilde{\nabla}U = (\nabla U, U_z) \) denotes the gradient in \( M \times (0, +\infty) \). This is a Hilbert space with the natural inner product that induces the norm above. Moreover, basically by definition, any \( U \in \tilde{H}^1(\tilde{M}) \) leaves a trace in \( L^2(M \times \{0\}) \).

\subsection{Applications of the extention method}
The nonlocal nature of the fractional Laplacian posses significant difficulties for its analysis. To overcome this challenge,
the {\em Caffarelli and Silvestre's extension method } transforms  nonlocal problems into  equivalent local ones in higher-dimensional space and has been applied successfully to study equations involving the fractional Laplacian, leading to a series of significant results (see \cite{BCPS} \cite{CZ}    and the references therein).

In \cite{BCPS}, among many interesting results, the authors investigated the properties of the positive solutions to
\begin{equation} (-\lap)^{s} u = u^p (x), \;\; x \in \mathbb{R}^n.
\label{e}
\end{equation}
They first used the {\em extension method} to transform the above nonlocal problem into a local one for $U(x,y)$ defined in the higher-dimensional half space $\mathbb{R}^n\times[0, \infty)$. Then by applying the {\em method of moving planes} on the local equation, they established the symmetry of $U(x,y)$ with respect to the variable $x$, and consequently proved  the non-existence of positive solutions in the subcritical case:
\begin{pro} \cite{BCPS}
Let $1/2 \leq s < 1$. Then the problem
\begin{equation}
\left\{\begin{array}{ll}
div(y^{1-2s} \nabla U)=0, & (x,y) \in \mathbb{R}^n\times[0, \infty),\\
\displaystyle-\lim_{y \ra 0^+} y^{1-2s} \frac{\partial U}{\partial y} = U^p (x,0), & x \in \mathbb{R}^n
\end{array}
\right.
\label{U}
\end{equation}
has no positive bounded solution provided $p < (n+2s)/(n-2s).$
\end{pro}

They then took trace to obtain
\begin{cor} Assume that $1/2 \leq s <1$ and $1<p < \frac{n+2s}{n-2s}$. Then equation (\ref{e}) possesses no bounded positive solution.
\label{mcor}
\end{cor}

A similar {\em extension method} was adapted in \cite{CZ} to obtain the nonexistence of positive solutions for an indefinite fractional problem:

\begin{thm} (Chen-Zhu)
Let $1/2 \leq s <1$ and $1<p<\infty$. Then the equation
\begin{equation}
(-\lap)^{s} u = x_1 u^p, \;\; x \in \mathbb{R}^n
\label{x1e}
\end{equation}
possesses no positive bounded solutions.
\end{thm}

The common restriction $s \geq 1/2$ is due to the approach that they need to employ  the {\em method of moving planes} on the solutions $U$ to the extended problem
\begin{equation}
div(y^{1-2s} \nabla U)=0, \;\; (x,y) \in \mathbb{R}^n\times[0, \infty).
\label{eeq}
\end{equation}
 Due to the presence of the weight factor $y^{1-2s}$, it is necessary to assume $s \geq 1/2$, and this condition appears to be indispensable for applying the {\em method of moving planes} to extended equation (\ref{eeq}). However, such a restriction does not arise in equation (\ref{x1e}), hence one may expect to remove the assumption when working directly with the original fractional equation. This is indeed the case as the readers will see when we introduce the {\em direct method of moving planes} in Section 5.
 \bigskip

\section{ Some maximum principles for elliptic equations}
\medskip

Maximum principles are among the most fundamental and powerful tools in the analysis of partial differential equations. They form the cornerstone for developing regularity theory, establishing Liouville-type theorems, and studying qualitative properties of solutions. In particular, the method of moving planes, the method of moving spheres, and the sliding method can all be viewed as continuous applications of suitable maximum principles.

\subsection{A simple maximum principle} \,
\smallskip

\begin{thm}(see \cite{Si} or \cite{CLL})
Assume that $\Omega$ is a bounded domain in $\mathbb{R}^n$, $u \in C^{1,1}_{loc} \cap \mathcal{L}_{2s}$ is lower semi-continuous on $\bar{\Omega}$ and satisfies
$$
\left\{\begin{array}{ll}
(-\lap)^{s} u(x) \geq 0 & x \in \Omega \\
u(x) \geq 0 & x \in \Omega^c.
\end{array}
\right.
$$
Then
\be \label{1.4}  u(x) \geq 0 \;\;\; x \in \Omega. \ee

We also have the following strong maximum principle:

Either $u(x)>0$ in $\Omega$ or $u(x)\equiv 0$ in the {\bf entire space} $\R^n$.

These conclusions hold for unbounded domain $\Omega$ if we further assume that
$$\underset{|x| \ra \infty}{\underline{\lim}} u(x) \geq 0.$$
\label{thm1.1}
\end{thm}

\begin{rem}
This {\em maximum principle} also holds if the fractional Laplacian $(-\lap)^s$ is replaced by the uniformly fractional operator $(-\lap)^s_a$, the nonlocal
Monge-Amp\`ere operator $D_s^\theta$, and the fractional p-Laplacian $(-\lap)^s_p$ mentioned in the Introduction.
\end{rem}

{\em Notice that, here $u(x^o)=0$ at one point in $\Omega$  implies $u(x) \equiv 0$ throughout the entire space $\R^n$. This is a key characteristic of nonlocal operators.
In contrast, if these nonlocal operators is replaced by the classical Laplacian $-\lap$ (a local operator), then the condition $u(x^o)=0$ has no influence at all on the values of $u$ outside $\Omega$. }
\bigskip

\subsection{A maximum principle for anti-symmetric functions} \,
\medskip

In applying the {\em method of moving planes}, what we consider are anti-symmetric functions. In this process, the corresponding {\em maximum principles} play important roles.

Let
$$T_\la =\{x \in \mathbb{R}^{n}\mid x_1=\lambda, \mbox{ for some } \lambda\in \mathbb{R}\}.$$
Let $$x^\la=(2\lambda-x_1, x_2, ..., x_n)$$
be the reflection of $x$ about the plane $T_\la$. Denote
$$\Sigma_\la =\{x \in \mathbb{R}^{n}\mid x_1<\lambda\},$$
and
$$u_\la(x) = u(x^\la), \;\; w_\la (x) = u_\la (x) - u(x).$$

\begin{thm} \label{mthm2} ( A maximum principle for anti-symmetric functions).

Let $\Omega$ be a bounded domain in $\Sigma_{\lambda}$.
Assume that $u \in \mathcal{L}_{sp}\cap C^{1,1}_{loc}$. If
\be
\left\{\begin{array}{ll}
(-\lap)^s_p u_{\lambda}(x) - (-\lap)^s_p u(x) \geq0  &\mbox{ in } \Omega,\\
w_{\lambda}(x) \geq0&\mbox{ in } \Sigma_{\lambda} \backslash\Omega
\end{array}
\right.
\label{s32}
\ee
with $p \geq 2$.
Then
\be
w_{\lambda}(x) \geq0 \mbox{ in } \Omega.\label{s41}
\ee

If $w_{\lambda} = 0$ at some point in $\Omega$, then
$$ w_\lambda(x) = 0 \; \mbox{almost everywhere in }  \mathbb{R}^n. $$

These conclusions hold for unbounded domain $\Omega$ if we further assume that
$$\underset{|x| \ra \infty}{\underline{\lim}} w_\lambda(x) \geq 0.$$
\end{thm}

\begin{rem}
Notice that here the exterior condition $w_\la(x) \geq 0$ is prescribed only in the complement of $\Omega$ in $\Sigma_\la$ instead of
the complement of $\Omega$ in the entire space, because $w_\la (x)$ is an anti-symmetric function.
\end{rem}
\medskip

\subsection{Maximum principle in a punctured ball and B\^{o}cher type theorem} \,
\medskip

In \cite{LWX}, the authors established maximum principles for fractional super harmonic functions in punctured balls and then connect these with
the fractional version of B\^{o}cher type theorems. These type of maximum principles are very useful in studying singular solutions, which appear naturally
in many physical and geometric phenomena. In particular, when we apply Kelvin transform, say $v(x) = \frac{1}{|x|^{n-2s}} u(\frac{x}{|x|})$, the transformed function
$v(x)$ develops a singularity at the origin.

\begin{thm} \label{LWX}  ( \cite{LWX} Fractional maximum principle in a punctured ball).

Suppose $n \geq 2$. Assume that $v \in \mathcal{L}_{2s}$ is a distributional nonnegative solution of
$$\left\{\begin{array}{ll}
(-\lap)^s v(x) + a(x) v(x) \geq 0, & x \in B_r(0) \setminus \{0\}, \\
v(x) \geq m > 0, & x \in B_r(0) \setminus B_{r/2}(0)
\end{array}
\right.
$$
with $a(x) \leq D$ for some constant $D$.

Then there exists a positive constant $c<1$ depending on $n$, $s$, and $D$ only, such that
$$ v(x) \geq c m, \;\; \forall \, x \in B_r(0) \setminus \{0\}, \forall \, r \leq 1. $$

\end{thm}

A similar maximum principle for anti-symmetric functions were also established in the same paper.
\medskip

\leftline{\bf An application to the method of moving planes.}
\smallskip

Consider
$$(-\lap)^s u(x) = u^p(x), \;\; x \in \R^n$$
with $1 <p \leq \frac{n+2s}{n-2s}$. In order to avoid imposing any decay condition on the solution $u$ near infinity,
one may apply the Kelvin transform
$$v(x) = \frac{1}{|x|^{n-2s}}u(\frac{x}{|x|^2}),$$
and then use {\em the method of moving planes} to study
$$w_\la (x) = v_\la(x) -v(x).$$

For $\la < 0$, the function $w_\la (x)$ has a singularity at $0^\la$, the reflection of the origin $0$ with respect to the plane $T_\la$ (see Figure 1).

{\begin{center}
\begin{tikzpicture}[scale=2.5]
    \def\lamval{1.0}  
    \def\r{0.6}       

    \draw[thick] (-3,0) -- (1.2,0);
    \draw[->, thick] (1.2,0) -- (1.5,0) node[below] {$x_1$};

    \fill[black] (0,0) circle (1pt);
    \node at (0,0) [below] {$0$};

    \draw[blue, thick] (-\lamval,-1) -- (-\lamval,1);
    \node[blue] at (-\lamval,1.2) {$T_\lambda$};

    \coordinate (Olambda) at (-2*\lamval,0);
    \fill[red] (Olambda) circle (1pt);
    \node[red] at (Olambda) [below] {$0^{\lambda}$};

    \draw[red, thick] (Olambda) circle (\r);
    \draw[red, dashed, ->] (Olambda) -- ++(30:\r) node[midway, above] {$r$};
    \node at (-\lamval, -1.2) {Figure 1. $0^{\lambda}$ is a singularity of $w_\lambda$};
\end{tikzpicture}
\end{center}}
During the application of the {\em method of moving planes}, it is often necessary to rule out the possibility that $w_\la$ attains its non-positive
minimum in a neighborhood of the singular point  $0^\la$. Before the establishment of Theorem \ref{LWX}, this difficulty was handled by using an equivalent integral formulation (see \cite{CLM}). However, not all pseudo-differential equations admit such integral representation. Theorem \ref{LWX} provides  a powerful tool to address this issue directly.

Indeed, in such a process, one already has
$$ w_\la (x) \geq m > 0, \;\; x \in B_{2r}(0^\la) \setminus B_r(0^\la).$$
Applying Theorem \ref{LWX}, one obtains
$$ w_\la(x) \geq c m, \;\; \forall \, x \in B_{2r}(0^\la) \setminus \{0^\la\}. $$
Consequently,  $w_\la$ cannot attain its non-positive minimum near the singularity $0^\la$. And if such a minimum were achieved elsewhere,
the standard argument would then lead to a contradiction.
\medskip

 In \cite{LWX},  the authors connected these maximum principles with the fractional version of the
B\^{o}cher-type theorem and provided a unified proof for both classical and fractional B\^{o}cher theorems.

\begin{thm}
Let $v(x) \in \mathcal{L}_{2s}$ with $n>2s$ be a nonnegative solution to
$$(-\lap)^s v(x) + c(x) v(x) = f(x) \geq 0, \; x \in B_1(x) \setminus \{0\}, $$
for some $f \in L^1_{loc} (B_1(0) \setminus \{0\})$ and $c(x)$ bounded from above.

Then there exists constant $a \geq 0$ such that

(i) $v(x), f(x) \in L^1_{loc} (B_1(0)),$

(ii) $(-\lap)^s v(x) + c(x) v(x) = f(x) + a \delta_0$ in $B_1(0)$.
\end{thm}

\bigskip

\subsection{A maximum principle in unbounded domains} \,
\medskip

In the earlier versions of the simple maximum principles, the analysis was restricted to bounded domains.
For unbounded domains $\Omega$, one typically assumes  that
$$\underset{|x| \ra \infty}{\underline{\lim}} u(x) \geq 0.$$
This additional condition essentially reduce it to bounded cases, since the minimum of $u$ would be attained at some interior point of $\Omega$, allowing one to derive a contradiction by the same argument used for bounded domains.

In what follows, we introduce a maximum principle that genuinely applied to unbounded domains without imposing any asymptotic condition on the solutions.

\begin{thm} (\cite{CW1} The maximum principle in unbounded domains) \label{unboundMP}

 Let $\Omega\subset\mathbb{R}^n$ be an open set, possibly unbounded and disconnected, and satisfying
\begin{equation}\label{big-omega^C}
\lim_{R\rightarrow +\infty}\frac{|B_R(x)\cap\Omega^c|}{|B_R(x)|} \geq c_0 >0
\end{equation}
for any $x\in \mathbb{R}^n$. Assume that the lower semi-continuous function $u\in C^{1,1}_{\rm loc}(\Omega) \cap \mathcal{L}_{2s}$ on $\overline{\Omega}$ is bounded from below, and satisfies
\begin{equation}\label{Umodel}
\left\{\begin{array}{r@{\ \ }c@{\ \ }ll}
&&\left(-\Delta\right)^{s}u(x)\geq0, & \ \ \mbox{at the points in} \,\, \Omega\,\, \mbox{where} \,\,u(x)<0\,, \\[0.05cm]
&&u(x) \geq 0, & \ \ x\in\Omega^c\,,
\end{array}\right.
\end{equation}
then
$u(x)\geq 0$ in $\Omega$\,.
\end{thm}

\begin{rem}
Roughly speaking, condition \eqref{big-omega^C} indicates that the ``size'' of the complement of $\Omega$ is ``not too small" as compared to the ``size'' of $\Omega$ in the sense of limit of the ratio. This ensures that the condition $u \geq 0$ in $\Omega^c$ exerts a noticeable influence on the values of $u$ inside $\Omega$.
\end{rem}
\bigskip

\section{Method of moving planes in integral forms}
\medskip

\subsection{An overview} \,
\medskip

 The {\em method of moving planes in integral forms} is another useful tool for studying qualitative properties of solutions to fractional equations. Since its introduction in \cite{CLO}, this method has been extensively applied by numerous researchers to solve a wide range of problems arising from nonlinear PDEs, particularly, those involving fractional equations and higher order equations (see \cite{CLO1, CLO2, FC, Ha, HLZ, Lei, MC, MCL, ZCCY} and the references therein).

 One major advantage of this method is that it works indiscriminately for equations of any order less than the dimension $n$. For instance, consider the fractional equation
 \be \label{B5.4} (-\lap)^s u(x) = f(x, u(x)), \;\; x \in \R^n. \ee
 The {\em direct method of moving planes} is typically applicable only for $0<s<1$. When $s>1$, one usually needs to decompose the higher-order equation into a system of lower-order equations and then apply the method to the system. Nonetheless, if one can establish the equivalence between pseudo-differential equation \eqref{B5.4} and the integral equation
 \be \label{B5.4a}
 u(x) = \int_{\R^n} \frac{f(y, u(y))}{|x-y|^{n-2s}} d y,
 \ee
 then the {\em method of moving planes in integral forms} works for all $s$ between $0$ and $n/2$. In this framework, the same argument applies seamlessly, regardless of whether $s<1$ or $s>1$.

 A primary limitation of this approach, however, is that it applies only to equations that can be reformulated as integral equations-such as those involving the classical Laplacian or the fractional Laplacian. The operator must be linear and possess an explicit Green's function. To the best of our knowledge, for most nonlocal operators introduced in the Introduction-such as general uniformly elliptic fractional operators, the nonlocal Monge-Amp\`ere operator, and fully nonlinear nonlocal operators (including the fractional $p$-Laplacian)-no corresponding integral representations are available. Consequently, this renders integral methods inapplicable to equations involving such operators.

To illustrate how the {\em method of moving planes in integral forms} can be employed to establish the radial symmetry and monotonicity of positive solutions, we consider the simple example
\be \label{A4.1}  (-\lap)^{s}u(x)=u^{p}(x), \;\; x\in \mathbb{R}^n. \ee
\medskip

\subsection{Turn the pseudo-differential equation into an integral equation} \,
\medskip

In order to apply this method, one first needs to show that the fractional equation is equivalent to an integral equation. The following Liouville type theorem plays a crucial role in this process.

\begin{thm} \label{thmLiouvilleE}
Assume that $u \in \mathcal{L}_{2s} \cap C_{loc}(\mathbb{R}^n)$ is $s$-harmonic in the sense of distribution, and $u$ is either bounded from below or from above, then $u \equiv C$.
\end{thm}

The earlier proof was given in \cite{BKN}, and more general version were established  in \cite{CDL} and \cite{Fa}.

Then based on this Liouville Theorem, we can prove

\begin{thm} \label{pro5.2} Assume $u\in \mathcal{L}_{2s}$ is a locally bounded nonnegative solution of
\be \label{1000}
(-\lap)^s u(x) = f(x),\ x\in \mathbb{R}^n
\ee
in the sense of distribution. Suppose
\be f(x) \geq 0, \;\; \mbox{ for } |x| \mbox{ sufficiently large.}
\ee

Then it is also a solution of
$$
u(x)= c + \int_{\mathbb{R}^n}\frac{C_{n, s }f(y)}{|x-y|^{n-2s}} dy.
$$
\end{thm}

Here we only provide the idea of proof. For more details, please see \cite{CGL}.

{\em Outline of proof}. Suppose that $u$ is a positive solution of (\ref{1000}). Since we do not impose any asymptotic behavior of  $u$ near infinity to guarantee the  convergence of the integral,  we first need to show that
$$
\int_{\mathbb{R}^n}\frac{f(y)}{|x-y|^{n-2s}}dy<\infty.
$$

Consider the function
$$v_R(x) = \int_{B_R(0)} \frac{C_{n,s}f(y)}{|x-y|^{n-2s}} d y.$$

Applying the {\em maximum principle}, we have, for $R$ sufficiently large,
$$ u(x) \geq v_R(x), \;\; x \in \R^n.$$

Letting $R \to \infty$, we arrive at, for each $x \in \R^n$,
$$ u(x) \geq  \int_{\R^n} \frac{C_{n,s}f(y)}{|x-y|^{n-2s}} d y \equiv v(x).$$

Now applying the Liouville Theorem (Proposition \ref{thmLiouvilleE}), to nonnegative function $w(x) = u(x) - v(x)$, we obtain
$$ u(x) = c + \int_{\R^n} \frac{C_{n,s}f(y)}{|x-y|^{n-2s}} d y,$$
for some nonnegative constant $c$.

Under appropriate conditions, for instance $f(y) = u^p(y)$, we can show that the constant $c$ must be $0$. Then we have
\begin{cor}
If $u$ is a nonnegative solution of \eqref{A4.1}, then it also satisfies the integral equation
$$ u(x) = \int_{\R} \frac{C_{n,s} u^p(y)}{|x-y|^{n-2s}} d y. $$
\end{cor}

\medskip

\subsection{Radial symmetry and monotonicity of solutions} \,
\medskip

Here we employ a simple example to illustrate how to obtain the radial symmetry and monotonicity of solutions by using {\em the method of moving planes in integral forms}.

Assume that $u$ is a nonnegative solution of the integral equation
\be \label{B5.2}
u(x)=\int_{\mathbb{R}^n}\frac{1}{|x-y|^{n-2s}}u^{p}(y)dy.
\ee
We will show that it is radially symmetric and monotone decreasing about some point in $\R^n$.

\medskip

\leftline{\em Outline of proof}
\smallskip

Let $T_\la$, $\hsp$, $x^\la$, and $w_\la (x)$ be defined as in the previous section. The general main frame is the same as that in the direct method of moving planes.

In step 1, we show that for $\la$ sufficiently negative
\be \label{B5.3} w_\la \geq 0, \;\; \forall \, x \in \hsp.
\ee
Previously, this was obtained by a {\em maximum principle}. Here for the integral equation, we approach differently by showing that the set
where inequality \eqref{B5.3} is violated must be empty, that is, we consider $$\hspn = \{ x \in \hsp \mid w_\la (x) < 0\}$$
where inequality \eqref{B5.3} is violated. We show that $\hspn$ is empty. To this end, we employ the HLS-inequality and the H\"{o}lder inequality to derive, for some $q>1$, an estimate of the form
\be \label{B5.4b}  \|w_\la\|_{L^q (\hspn)} \leq (\cdots) \|w_\la\|_{L^q (\hspn)}. \ee

For $\la$ sufficiently negative, under suitable integrability condition, we show that the quantity in $(\cdots)$ is small (actually, $<1$ is sufficient to serve the purpose), hence
$\|w_\la\|_{L^q (\hspn)} = 0$. This verifies \eqref{B5.3}.

Same estimate \eqref{B5.4b} works in the second step with $\hspn$ being a subset containing in a neighborhood of infinity and a narrow region.

Here we skip the details. Interested readers may see Chapter 7 of the book \cite{CLM}.
\medskip

In \cite{ZCCY}, the authors consider nonnegative solutions to a system of $m$ fractional equations
\be
(-\lap)^s u_i(x) = f_i(u_1(x), \cdots u_m(x)), \;\; i = 1, \cdots, m.
\label{A500}
\ee
They first established the equivalence between the pseudo-differential system \eqref{A500} and the integral system
\be
u_i(x) = \int_{\R^n} \frac{c_{n,s}}{|x-y|^{n-2s}} f_i(u_1(y), \cdots u_m(y)) d y,  \;\; i = 1, \cdots, m.
\label{A501}
\ee
Then by using the {\em method of moving planes in integral forms}, they obtained radial symmetry of solutions in the critical case and nonexistence in the
subcritical case.
\bigskip

\section{The direct methods}
\medskip

As mentioned in the Introduction, both the extension method and the integral equation method were primarily developed for equations involving the fractional Laplacian. Before applying these approaches, one must either extend the problem to a higher-dimensional local formulation or derive an equivalent integral equation. In these processes, certain additional conditions are often required. Moreover, to the best of our knowledge, these methods are generally not applicable to equations involving other types of nonlocal operators, such as uniformly elliptic fractional operators, fractional Monge-Amp\`ere operators, and fractional $p$-Laplacians.

In recent years, several direct methods have been introduced, including the direct methods of moving planes, moving spheres, blow-up and rescaling, and the sliding method. In this section, we will use simple examples to illustrate the essential ideas behind these approaches.

\subsection{Methods of moving planes} \,
\medskip

Here, we use a simple example to explain how to apply the method of moving planes to obtain radial symmetry and monotonicity of positive solutions for a fractional equation.

\subsubsection{An example} \,
\smallskip

Consider
\begin{equation}
\left\{ \begin{array}{ll}
(-\lap)^{s} u(x)  = f(u(x)), u(x) >0,  & x \in B_1(0), \\
u(x) = 0 , & x \in B_1^c(0).
\end{array}
\right.
\label{s65}
\end{equation}
We will show that every solution $u$ is radially symmetric and monotone decreasing about the origin. More precisely, we have
\begin{thm} \label{thm4.1}
Assume that $u \in C_{loc}^{1,1}(B_1(0)) \cap C(\overline{B_1(0)})$ is a positive solution of (\ref{s65}) with $f(\cdot)$ being Lipschitz continuous. Then $u$ must be radially symmetric and monotone decreasing about the origin.
\end{thm}

\smallskip

\leftline{\em Outline of the approach}
\smallskip

 Place $x_1$ axis in any given direction. Let $T_\la$, $\Sigma_\la$, and $x^\la$ as previously defined (see Figure 2).

 {\begin{center}
\begin{tikzpicture}[scale=2]
    \def\lamval{0.3}    
    \def\xpos{-0.6}     
    \def\sphereR{1}   

    \draw[thick] (-2,0) -- (2,0);
    \draw[->, thick] (2,0) -- (2.2,0) node[below] {$x_1$};

    \fill[black] (0,0) circle (1pt);
    \node at (0,0) [below] {$0$};

    \draw[blue, thick] (-\lamval,-1.2) -- (-\lamval,1.2);
    \node[blue] at (-\lamval,1.4) {$T_\lambda$};

    \node at (-\lamval-1, 0.8) {$\Sigma_\lambda$};

    \coordinate (x) at (\xpos,0.5);
    \coordinate (xlambda) at (-\lamval*2-\xpos,0.5);

    \fill[red] (x) circle (1pt);
    \node[red] at (x) [above] {$x$};

    \fill[red] (xlambda) circle (1pt);
    \node[red] at (xlambda) [above] {$x^{\lambda}$};

    \draw[dashed] (x) -- (xlambda);

    \draw[thick, orange] (0,0) circle (\sphereR);

    \node at (-\lamval, -1.6) {Figure 2.};
\end{tikzpicture}
\end{center}}

 To derive the radial symmetry of solution $u(x)$ about the origin, since the direction $x_1$ is chosen arbitrarily, it suffice to show that it is symmetric about the plane $T_0$ passing through the origin. This is equivalent to proving that $w_0(x) \equiv 0$.

To achieve this goal, we first place the plane $T_\la$ very close to the left end of the ball, that is, $\lambda$ is close to $-1$. In this position,
$\Sigma_\la$ is a narrow region, allowing us to apply the {\em narrow region principle} to conclude that $w_\la (x) \geq 0$. This provides a starting point for moving the plane. There has been a similar {\em narrow region principle} for local equations. It was proved by constructing an auxiliary function $\bar{w} = \frac{w_\la}{g}$, which no longer work in the fractional case. For the fractional equation, we apply {\em 2s-order derivative like} property (see \cite{CGL}) of $w_\la (x)$ at its minimum point which is derived from the singular integral representing the fractional Laplacian and the anti-symmetry property of $w_\la(x)$.

We then move the plane $T_\la$ continuous to the right
as long as the inequality holds, until it reaches its limiting position $T_{\bar{\la}}$. We show that $\bar{\la}$ must be $0$ by using a contradiction argument.

At each moment, the plane is moved forward slightly, and $\Sigma_\la$ is increased by a narrow region. Thus, this process of moving planes can be interpreted as applying the maximum principle (the narrow region principle) continuously for infinitely many times. For more details, please see \cite{CGL}.
\medskip

\leftline{\em A key ingredient}
\smallskip

Since in the process of moving planes, what we considered are anti-symmetric functions, we need the following {\em narrow region principle for antisymmetric functions.}

\begin{lem}

Let $\Omega$ be a bounded narrow region in
$$\Sigma_\lambda \equiv \{ x \in \mathbb{R}^n \mid x_1 < \lambda \}, $$ such that it is contained in the narrow slab
 $$\{x \mid  \lambda-\delta<x_1<\lambda \, \}$$
 for some $\delta > 0$. Suppose  that $w_\lambda \in {\mathcal L}_{2s} \cap C_{loc}^{1,1}(\Omega)$ and is lower semi-continuous on $\bar{\Omega}$. If
 $c(x)$ is bounded from below in $\Omega$  and
\be
\left\{\begin{array}{ll}
(-\lap)^s w_\lambda (x) +c(x)w_\lambda(x)\geq0  &\mbox{ in } \Omega,\\
w_\lambda(x) \geq0&\mbox{ in }  \Sigma_\lambda \backslash\Omega,
\end{array}
\right.
\label{wcg0}
\ee
then for sufficiently small $\delta$, we have
$$
w_\lambda (x) \geq0 \mbox{ in } \Omega.
$$
These conclusions hold for unbounded region $\Omega$ if we further assume that
$$\underset{|x| \ra \infty}{\underline{\lim}} w_\lambda(x) \geq 0 .$$

The following strong maximum principle is valid for any bounded domain $\Omega$ (need not to be narrow):

If it is shown that
$$
w_\lambda (x) \geq0 \mbox{ in } \Omega,
$$
then we have
$$\mbox{ either } w_\la (x) > 0 \; \mbox{ in } \Omega \; \mbox{ or } w_\la (x) \equiv 0 \; \mbox{ in } \R^n.$$

\label{lem4.1}
\end{lem}

\begin{rem}  This narrow region principle may be view as a ``better'' maximum principle in the following sense. As one have seen at the beginning of Section 4, the maximum principle holds if
$$ (-\lap)^s w_\lambda (x) \geq0  \mbox{ in } \Omega.$$
We may think $(-\lap)^s$ as a ``positive'' operator. Now if we add a nonnegative function $c(x)$, then $(-\lap)^s + c(x)$ is still a ``positive'' operator, and the maximum principle is also valid. However, in practice we do not know the sign of $c(x)$ and hence cannot establish the maximum principle in a bounded region $\Omega$ for the operator $(-\lap)^s + c(x)$. Here comes the remedy: to make the region narrow. For such a narrow region, we only require $c(x)$ to be bounded from below to ensure the maximum principle.
\end{rem}

In \cite{WN}, based on the {\em direct method of moving planes}, and by a novel approach, the authors established a  similar {\em narrow region principle} for the fractional $p$-Laplacian.

\begin{lem}

Let $\Omega$ be a bounded narrow region in
$$\Sigma_\lambda \equiv \{ x \in \mathbb{R}^n \mid x_1 < \lambda \}, $$ such that it is contained in the narrow slab
 $$\{x \mid  \lambda-\delta<x_1<\lambda \, \}$$
 for some $\delta > 0$. Suppose  that $w_\lambda \in {\mathcal L}_{2s} \cap C_{loc}^{1,1}(\Omega)$ and is lower semi-continuous on $\bar{\Omega}$. If
 $c(x)$ is bounded from below in $\Omega$  and
\be
\left\{\begin{array}{ll}
(-\lap)^s_p  u_\lambda (x)  - (-\lap)^s_p u(x) +c(x)w_\lambda(x)\geq0  &\mbox{ in } \Omega,\\
w_\lambda(x) \geq0&\mbox{ in }  \Sigma_\lambda \backslash\Omega,
\end{array}
\right.
\ee
and there exists $y^0 \in \Sigma_\la$ such that $w_\la(y^0)>0$.

Then for sufficiently small $\delta$, we have
$$
w_\lambda (x) \geq0 \mbox{ in } \Omega.
$$
These conclusions hold for unbounded region $\Omega$ if we further assume that
$$\underset{|x| \ra \infty}{\underline{\lim}} w_\lambda(x) \geq 0 .$$

The following strong maximum principle is valid for any bounded domain $\Omega$ (need not to be narrow):

If it is shown that
$$
w_\lambda (x) \geq0 \mbox{ in } \Omega,
$$
then we have
$$\mbox{ either } w_\la (x) > 0 \; \mbox{ in } \Omega \; \mbox{ or } w_\la (x) \equiv 0 \; \mbox{ in } \R^n.$$
\end{lem}

Employing this {\em narrow region principle}, the authors in \cite{WN} obtained the same radial symmetry for the solution of \eqref{s65} with the fractional
Laplacian $(-\lap)^s$ replaced by the fractional $p$-Laplacian $(-\lap)^s_p$.

For the equation involving the fractional $p$-Laplacian in the whole space
\be
(-\lap)^s_p u(x) = g(u(x)), \;\; x \in \mathbb{R}^n .
\label{meqws}
\ee

The authors in \cite{CL2} and \cite{WN} also established radial symmetry and monotonicity via {\em direct method of moving planes}. One of them reads as
\begin{thm}
Assume that $u \in C^{1,1}_{loc} \cap \mathcal{L}_{sp}$ is a positive solution of (\ref{meqws}) with
$\lim_{|x| \ra \infty} u(x) = 0.$
Suppose $g'(s) \leq 0$ for $s$ sufficiently small.

Then $u$ must be radially symmetric and monotone decreasing about some point in $\mathbb{R}^n$.
\label{mthm5}
\end{thm}
\medskip

\subsection{Method of moving spheres} \,
\medskip

\subsubsection{Outline of the method} \,
\smallskip

 In the {\em method of moving planes}, we move a family of parallel planes $T_{\lambda}$ along a chosen direction toward a limiting position to obtain symmetry of the solutions with respect to the limiting plane. In contrast, while in the {\em method of moving spheres}, we fix a center of the sphere and
 compare the values of a solution $u$ with its Kelvin transform.  Then we gradually increase or decrease the radius of the spheres to derive
 some kind of monotonicity--specifically, the comparison between $u$ and its Kelvin transform--along the radial directions of the spheres. If such monotonicity holds for arbitrary centers, one can deduce non-existence of solutions; or if such symmetry holds for arbitrary centers, all solutions can be classified. However, this approach does not yield the usual notions of monotonicity and symmetry of solutions.

 An advantage of the method of moving spheres is that, when classification of solutions is possible, it often yields the result directly. By contrast, when applying the method of moving planes, one must first establish the radial symmetry of the solutions and then employ additional arguments to determine their explicit forms.

Since the method of moving spheres relies essentially on the Kelvin transform, it is mainly applicable to equations involving the classical Laplacian or the fractional Laplacian. For other types of operators-such as those discussed in Section 3, including general uniformly elliptic operators, nonlocal Monge-Amp\`ere operators, and fully nonlinear nonlocal operators such as the fractional $p$-Laplacian-the method becomes ineffective.

 To illustrate the general idea, we consider a  simple example below.

 \subsubsection{Classifications of solutions for a Lane-Emden equation}  \,
 \medskip

In \cite{CLZ} and in \cite{LXYZ}, applying the {\em direct method of moving spheres}, the authors obtained the following
 \begin{pro}\label{classification-theorem0}
Assume that $g:\dR_+^1\rightarrow \dR_{+}^1\cup\{0\}$
is locally bounded and
$\frac{g(r)}{r^{p}}$ is non-increasing with $p = \frac{n+2s}{n-2s}$.
If $u\in \mathcal{L}_{\alpha}\cap C_{loc}^{1,1}(\dR^n)$ is a nonnegative solution to equation
$$(-\triangle)^{s}u(x)=g(u(x)), \ x\in\dR^n, $$
then one of the following holds:

$(1)$ For some $C_0\geq0$, $u(x)\equiv C_0$ for every $x\in\dR^n$ and $g(C_0)=0$.

$(2)$ There exists $\gamma, \la >0$, $x_0\in\dR^n$ such that
\begin{equation}
u(x)=\frac{\gamma}{(|x-x_0|^2+\la^2)^{\frac{n-2s}{2}}},\ \forall x\in\dR^n,
\end{equation}
and $g(r)$ is a multiple of $r^{p}$
for every $r\in(0,\max\limits_{x\in\dR^n}u(x)]$.
\end{pro}

\begin{rem} Although in the statement of the theorem in \cite{CLZ}, the function $g(\cdot)$ is not required to be monotone increasing, one can see from their key estimate (\cite{CLZ}, Lemma A.2), this condition is indeed needed. In fact, in the proof of this Lemma A.2, the authors relied on an integral equation and take advantages of monotone increasing property of $t^p$.

Later, in \cite{LXYZ}, the authors removed this monotonicity condition by employing two new maximum principles established there (see Theorems 1.1 and 1.2 in \cite{LXYZ}). Moreover, their proof of the theorem
does not depend at all on the integral representation of the solution, and hence this method can
also be used for more general non-local equations, where no equivalent integral equations exist.
\end{rem}

In the special case of the Lane-Emden type equation
\be (-\lap)^{s} u(x) = u^p(x), \; \; x \in \mathbb{R}^n,
\label{LEE}
\ee
the above theorem implies immediately

\begin{cor}

i) In the subcritical case $1<p<\frac{n+2s}{n-2s}$, equation \eqref{LEE} does not admit positive solutions.

ii) In the critical case $p=\frac{n+2s}{n-2s}$, all the positive solutions must assume the form
$$ u(x)=\frac{\gamma}{(|x-x_0|^2+\la^2)^{\frac{n-2s}{2}}},\ \forall x\in\dR^n.$$

\end{cor}

 One can use the {\em method of moving spheres} to prove the nonexistence and classification of positive solutions as outlined in the following.
\medskip

 Taking a point  $x_0\in\dR^n$ as the center and any real number $\lambda>0$ as the radius, we define the Kelvin transform of
 $u$ centered at $x_0$ with respect to a sphere of radius $\lambda$ as follows,
 \begin{equation}
 u_{\lambda}(x)\equiv(\frac{\lambda}{|x-x_0|})^{n-2s}u(x^{\lambda}),
\label{Kelvin-transform}
\end{equation}
where $$x^{\lambda}\equiv\frac{\lambda^2 (x-x_0)}{|x-x_0|^2}+x_0$$
is the inversion of $x$ with respect to the sphere
$$S_{\lambda}(x_0) \equiv \{ x \mid |x-x_0| = \lambda \}.$$

For simplicity of writing, here we take $x_0=0$.

We compare the value of $u(x)$ with its Kelvin transform $u_{\lambda}(x)$ with respect to the sphere
$S_{\lambda}(0)$ centered at the origin. Notice that for $x \in B_{\lambda}(0)$, $u_{\lambda}(x)$ actually is the value of $u$ at $x^{\lambda}$ outside the ball multiplied by $\left(\frac{\lambda}{|x|}\right)^{n-2s}$.

Let $w_\la(x)= u_{\lambda}(x) - u(x)$, then one can verify that $w_\la$ is anti-symmetric in the sense
$$ \left( \frac{\lambda}{|x|}\right)^{n-2s}w_{\lambda}(x^\lambda) =-w_{\lambda}(x).$$
One can also show that $w_\la$ satisfies the equation
$$(-\triangle)^{s}w_\la(x)+c_\la(x)w_\la(x)= 0,$$
for some $c_\la(x)$ depending on $u(x)$.

In step 1, we show that for $\lambda$ sufficiently small, it holds
\be w_\la(x) \geq 0 , \;\; \forall \, x \in B_{\lambda}(0).
\label{8.1a}
\ee
This will be accomplished by the {\em narrow region principle:}

\begin{lem} \label{NRS}
Let $w_\la \in \mathcal{L}_{2s}\cap C_{loc}^{1,1}(\Omega)$ be lower
semi-continuous on $\bar{\Omega}$. If $c(x)$ is bounded from below
in $\Omega$ and
\begin{equation*}
\left\{\begin{array}{ll}
(-\triangle)^{s}w_\la(x)+c(x)w_\la(x)\geq 0,& x\in\Omega\subset B_{\lambda}(0),\nonumber\\
w_\la(x)\geq 0,& x\in B_{\lambda}(0)\setminus\Omega,
\nonumber\\
\left( \frac{\lambda}{|x|}\right)^{n-2s}w_\la(x^\la)=- w_{\lambda}(x),&x\in B_{\lambda}(0).
\end{array}
\right.
\end{equation*}
 If $\Omega$ is contained in a
spherically narrow region (see Figure 3), i.e.,
$$\Omega\subset A_{\lambda-\delta,\lambda}(0) \equiv
\{x\in\dR^n \mid \lambda-\delta<|x|<\lambda\},$$
for sufficiently small $\delta>0$,
then we have
\begin{equation*}
\inf\limits_{x\in\Omega}w_\la(x)\geq 0.
\end{equation*}
Furthermore, if $w_\la(x)=0$ for some $x\in\Omega$, then $w_\la(x)\equiv0$ in
$\dR^n$.
\end{lem}

\begin{center}
\begin{tikzpicture}[scale=2]
    \def\R{1}        
    \def\r{0.8}      
    \def\xangle{30}  
    \def\xdist{0.5}  


    \draw[thick] (0,0) circle (\R);
    \draw[thick] (0,0) circle (\r);

    \fill[gray!30, opacity=0.5] (0,0) circle (\R);
    \fill[white] (0,0) circle (\r);

    \fill[black] (0,0) circle (1pt);
    \node at (0,0) [below left] {$0$};

    \coordinate (x) at (\xangle:\xdist);
    \coordinate (xlambda) at (\xangle:\R*\R/\xdist); 

    \fill[red] (x) circle (1pt);
    \node[red] at (x) [above right] {$x$};

    \fill[blue] (xlambda) circle (1pt);
    \node[blue] at (xlambda) [above right] {$x^{\lambda}$};

    \draw[dashed, red] (x) -- (0,0) -- (xlambda);

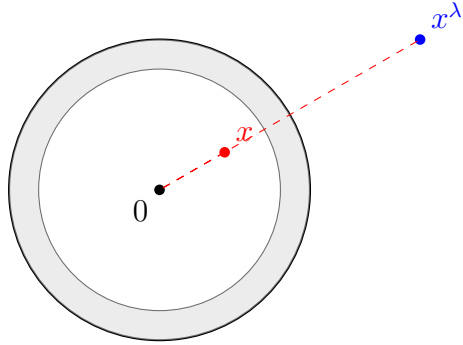
    \node at (0, -1.6) {Figure 3. Spherical narrow region.};

\end{tikzpicture}

\end{center}

For sufficiently small $\la$, $B_{\lambda}(0)$ is a spherical narrow region, and hence the above lemma implies \eqref{8.1a}.

Step 1 provides a starting point to continuously dilate the sphere with the same center. Now in Step 2,  we increase $\lambda$, the radius of the sphere,
as long as inequality (\ref{8.1a}) holds. Define
$$ \lambda_o = \sup\{ \lambda \mid w_{\mu}(x) \geq 0, \, x \in B_{\mu}(0); \forall \, \mu \leq \lambda \}.$$

(a) {\em In the subcritical case,} we can show that $\lambda_o = \infty,$ by employing the {\em narrow region principle} again. Because each time we only increase the radius $\lambda$ a little bit, the region in between the two spheres is a spherically narrow region.
These steps can be carried through for any center $x_0$, hence we have actually proved that
$$(\frac{\lambda}{|x-x_0|})^{n-2s}u(x^{\lambda}) \geq u(x) , \;\; \forall \, \lambda >0.$$
This is a kind of the monotonicity we mentioned earlier.

Due to the arbitrariness of the center $x_0$ and the radius $\la$, by a remarkable general result of Li and Zhu \cite{LZh}, one will be able to deduce that
$$\nabla u(x) = 0 , \quad \forall x \in \mathbb{R}^n.$$
Consequently, $u$ must be constant, which contradicts the equation
 if $u$ is a positive solution.

(b) {\em In the critical case}, a similar argument as in the above will lead to $\la_0 < \infty$ for every given center. Then again by a result in \cite{LZh}, one will arrive at the classification of all positive solutions.

This sketches the proof. A detailed argument can be found in \cite{CLZ} or in Section 8.3 of \cite{CLM}.
\medskip

\subsubsection{Method of scaling spheres}

In \cite{DQ1}, another useful method was introduced to prove the nonexistence solutions--a Liouville type theorem. It is called
the {\em method of scaling spheres}.

Being distinct from the {\em method of moving spheres}, the {\em method of scaling spheres} fully exploits the integral representation of solutions together with the {\em bootstrap} technique to derive sharp lower bound estimates for the asymptotic behavior of positive solutions near singularities or at infinity. This approach ultimately yields Liouville-type theorems. In addition to improving most known Liouville-type results to their optimal exponents in a unified and simple manner, the method of scaling spheres can be conveniently applied to a wide range of problems, including those involving singularities or lacking translation invariance in general domains.

The method proceeds in three main steps.
First, the integral representation implies that any positive solution $u$ satisfies the same lower bound estimate as that of the fundamental solutions as
$|x| \to \infty$.

Second, by combining this lower bound with the integral representation and applying the scaling spheres procedure, one obtains an improved lower bound for $u$ stronger than that of the fundamental solution.

Finally, through a {\em bootstrap} iteration process, in conjunction with the integral representation, increasingly sharper lower bound estimates are established for the asymptotic behavior of $u$. This iterative improvement ultimately contradicts the integrability condition implied by the integral representation, unless the solution is identically zero.
\medskip

\subsection{The sliding method} \,
\medskip

The {\em sliding method} is mainly used to obtain monotonicity and one dimensional symmetry of solutions. For certain regions, such as half-spaces. it can also be utilized to derive the uniqueness of solutions.

\subsubsection{Different requirements on operators} \,
\smallskip

When applying the {\em methods of moving planes} to equations involving non-local operators defined by singular integrals
$$ \mathcal{L} u(x) = \int_{\R^n} k(x-y)[u(x)-u(y)] d y, $$
one requires the kernel of the operator $k(x-y)$ be monotone decreasing in $|x-y|$. This condition is obviously satisfied by the fractional Laplacian with kernel
$$k(x-y) = \frac{C}{|x-y|^{n+2s}}. $$

Nevertheless, for other nonlocal operators, such as the uniformly elliptic fractional operators and the nonlocal Monge-Amp\`ere operators,
this requirement is not met.
Consequently, the method of moving planes can no longer be applied to establish the qualitative properties of solutions to the nonlocal problem involving such operators. Nonetheless, it can often be realized by the sliding method. The key ingredient for implementing the sliding method is the establishment of the maximum principles for the solutions themselves rather than for anti-symmetric functions. As a result, this method can be applied to a number of different operators provided they obey suitable maximum principles.

However, this does not imply that the sliding method is superior to the method of moving planes. In fact, the sliding method usually fails to prove the symmetry of solutions and typically requires stronger restrictions on the nonhomogeneous terms as well as on the exterior or asymptotic conditions.
Despite such limitations, the sliding method remains a valuable tool in many applications.
\medskip

\subsubsection{A fractional version of Gibbon's conjecture} \,
\smallskip

As an application, we employ the {\em sliding method} to establish a fractional version of Gibbon's conjecture, a weaker version of
the well-known {\em De Giorgi's conjecture}, which states that
\smallskip

{\em Suppose that $u(x)$ is an entire solution of the equation
\begin{equation}\label{AC-equ}
  -\Delta u=u-u^3,\,\, x=(x',x_n)\in \mathbb{R}^n,
\end{equation}
satisfying
\begin{equation*}
  |u(x)|\leq 1\,\,\mbox{and}\,\, \frac{\partial u}{\partial x_n}>0 \,\, \mbox{in} \,\,\mathbb{R}^n.
\end{equation*}
Then the level sets of $u(x)$ must be hyperplanes, at least for $n\leq 8$.}
\smallskip

There have seen many partial results concerning this conjecture, so far it has not been validated completely.

Gibbons proposed a weaker version of De Giorgi's conjecture that replaces the one-direction monotonic condition $\frac{\partial u}{\partial x_n}>0$ by a uniform convergence assumption
\begin{equation}\label{uniform-convergence}
 \displaystyle\lim_{x_n\rightarrow\pm\infty}u(x',x_n) = \pm1 \,\,\mbox{uniformly with respect to} \,\, x'\in \mathbb{R}^{n-1}.
\end{equation}
This is referred to as

\textbf{Gibbons' conjecture.} Let $u(x)$ be an entire solution of \eqref{AC-equ} satisfying
$|u(x)|\leq 1$ and \eqref{uniform-convergence}, then $u(x)$ is monotonically increasing with respect to $x_n$ and depends exclusively on $x_n$.

Fortunately, Gibbons' conjecture has been successfully established in any dimensions by various methods, including the moving plane method by Farina \cite{Far}, the probabilistic arguments by Barlow, Bass, and Gui \cite{BBG}, and the sliding method by Berestycki, Hamel, and Monneau \cite{BHM}. It is worth mentioning that their results applied to more general nonhomogeneous terms $f$ than the De Giorgi-type nonlinearities $f(u)=u-u^3$.
Among them, the celebrated sliding method was first introduced by Berestycki and Nirenberg \cite{BN3} to establish qualitative properties of positive solutions to the local elliptic equations.

Afterwards, Wu and Chen \cite{WuC, WC2} developed a direct sliding method for fractional equations to establish monotonicity, one-dimensional symmetry, uniqueness, and nonexistence of solutions. Such direct method avoids the heavy use of classical extension method  introduced in \cite{CS} to overcome the difficulties caused by the non-locality of fractional operators. More importantly, this direct sliding method can be used to extend and prove Gibbons' conjecture in the settings of other fractional elliptic equations involving
various nonlocal operators.

\begin{thm} \cite{WuC} \label{thm3.2}
Let $u\in C^{1,1}_{\rm loc}(\mathbb{R}^n) \cap \mathcal{L}_{2s}$ be a solution of the problem
\begin{equation*}
\left\{\begin{array}{ll}
(-\Delta)^s u(x) = f(u(x)),\, |u(x)|\leq1,\,\, x \in \mathbb{R}^n,\\
\displaystyle\lim_{x_n\rightarrow\pm\infty}u(x',x_n) = \pm1, \,\mbox{uniformly for}\, x'=(x_1,...,x_{n-1})\in \mathbb{R}^{n-1}.
\end{array}
\right.
\end{equation*}
Assume that $f\in C([-1,1])$ is
 non-increasing on $[-1,-1+\delta]\cup [1-\delta,1]$ for some $\delta>0$.

 Then $u(x)$ is strictly increasing with respect to $x_n$ and depends only on $x_n$, i.e. $u(x)=u(x_n)$ for every $x=(x',x_n)\in\mathbb{R}^n$.
\end{thm}

\begin{rem}
It can be easily verified that if $f(u) = u-u^3$, then the conditions in the theorem are satisfied. Hence this theorem establishes a fractional version of Gibbons' conjecture.
\end{rem}

\begin{rem} If the fractional Laplacian $(-\lap)^s$ is replaced by the fractional p-Laplacian $(-\lap)^s_p$, the above theorem is still valid (see \cite{WC2}).
\end{rem}

\begin{proof}[The outline of proof of Theorem \ref{thm3.2}]

We use the sliding method.

Let  $x^\tau=x+\tau e_n$ with $e_n = (0, \cdots, 0, 1)$ and $u^\tau (x) = u(x^\tau)$. We compare the value of $u(x)$ with $u^\tau (x)$.
Intuitively, the graph of $u^\tau(x)$ is obtained by sliding the graph of $u(x)$ downward (in $-x_n$ direction) $\tau$ units. Due to the asymptotic behavior of $u$,
$$\displaystyle\lim_{x_n\rightarrow\pm\infty}u(x',x_n) = \pm1, \,\mbox{uniformly for}\, x'=(x_1,...,x_{n-1})\in \mathbb{R}^{n-1},$$
 one may expect that for sufficiently large $\tau$, it holds $u^\tau (x) \geq u(x)$. This will be verified by a {\em maximum principle}.
Then we continuously decrease $\tau$ (slide the domain upwards) as long as the inequality holds to its limiting position and show that such process can be continued until $\tau = 0$. This arrives at the monotonicity of $u$ in $x_n$-direction.

For convenience of applying the {\em maximum principle}, we consider the difference
$$w^\tau(x)= u^\tau(x)-u(x).$$

The proof is divided into three steps.
\medskip

\noindent \textbf{Step 1.} For sufficiently large $\tau$, we show that
\begin{equation}\label{4.2.1}
w^\tau(x)\geq0, \, x\in \mathbb{R}^n.
\end{equation}

This is done by applying the {\em maximum principle in unbounded domains} (Theorem \ref{unboundMP}) to the function
$$\bar{w}^\tau(x)=w^\tau(x)- \frac{1}{2} \inf_{x\in\mathbb{R}^n}w^\tau(x)$$
in the domain
$$ \Omega = \R^{n-1} \times (-\infty, M).$$
Here $M$ is chosen to be sufficiently large, so that
$$ \bar{w}^\tau(x) \geq 0, \;\; x \in \Omega^c.$$

Inequality \eqref{4.2.1} provides a starting point for sliding the domain along the $x_n$-axis.
\medskip

\noindent \textbf{Step 2.} We continuously decrease $\tau$ to its limiting position as long as inequality \eqref{4.2.1} holds and define
\begin{equation*}
  \tau_0:=\inf\{\tau\mid w^\tau(x)\geq0, \, x\in \mathbb{R}^n \}.
\end{equation*}
We claim $\tau_0=0$ by a contradiction argument. More precisely, if $\tau_0>0$,
we prove that $\tau_0$ can be decreased a little bit
while the inequality \eqref{4.2.1} is still valid, which contradicts the definition of $\tau_0$. Furthermore, by virtue of the strong maximum principle, we deduce that
\begin{equation*}
  w^\tau(x)>0, \, x\in \mathbb{R}^n \,\, \mbox{for any}
\, \tau>0.
\end{equation*}
Thus, we conclude that $u$ is strictly increasing
with respect to $x_n$.

\noindent \textbf{Step 3.} We finally prove that $u(x)=u(x',x_n)$ is independent of $x'$.
By proceeding similarly as in Step 1 and Step 2, we can derive
$$u(x+\tau\nu)>u(x),\  x\in \mathbb{R}^n, \, \forall\tau>0,$$
for every vector $\nu=(\nu_1,...,\nu_n)$ with $\nu_n>0$. This implies that $u(x)$ is strictly increasing along any direction which has an acute angle with the positive $x_n$ axis. Let $\nu_n \rightarrow 0$, then by the continuity of $u(x)$, the inequality is still preserved in the sense
$$u(x+\tau\nu)\geq u(x),\  x\in \mathbb{R}^n, \, \forall\tau>0.$$
Now $\nu$ can be any given direction perpendicular to $x_n$ axis. This implies that $u(x)$ is independent of $x'$. That is $u(x)=u(x_n)$.
\end{proof}

\medskip

\subsubsection{Comparison between the method of moving planes and the sliding method} \,
\smallskip

The method of moving planes is a powerful tool and has been employed extensively to study the qualitative properties of solutions to partial differential equations, including symmetry, monotonicity, and nonexistence. It can also be used to derive a priori bounds of solutions in a boundary layer of a domain.

In contrast, sliding methods,  while less popular, are primarily employed to derive monotonicity of solutions and, occasionally,  uniqueness.

When applying the method of moving planes to equations involving non-local operators defined by singular integrals, it is crucial that the operator's kernel satisfies certain monotonicity properties--specifically being monotone decreasing in $|x-y|$. This condition is met by the fractional Laplacian but not by several other operators, such as all the four operators listed in Section 3,
hence the method of moving planes is no longer applicable in such cases. In these scenarios, the sliding method becomes useful. Its key requirement is the establishment of the maximum principles for the solutions rather than for the anti-symmetric functions, eliminating the need for kernel monotonicity.  This makes the sliding method applicable to a broader class of operators, provided they obey maximum principles.

 However, this does not mean that the sliding method is better than the method of moving planes, indeed, the sliding method usually fails to prove the symmetry of solutions, and it needs to impose stronger restrictions on the nonhomogeneous terms and exterior or asymptotic  conditions.
Despite such limitations, the sliding method remains valuable for various  applications. Interested readers may see \cite{CHM} for more deliberations and examples.
\medskip

\subsection{Blow-up and rescaling} \,
\medskip

We will use the following simple example to demonstrate how the blow-up and rescaling analysis can be applied to establish a priori estimate for positive solutions.

Let $\Omega$ be a bounded domain in $\R^n$ with smooth boundary. Consider

\be\label{Aeq1.7}
\left\{\begin{array}{ll}
 (-\Delta)^s u(x)= f(x, u(x)), & x\in \Omega,\\
 u(x) \equiv 0, & x \in \Omega^c.
 \end{array}
 \right.
\ee

Assume that
\smallskip

 (f1) $f(x, t): \Omega \times [0, \infty) \to \mathbb{R}$ is uniformly H\"{o}lder continuous with respect to $x$ and continuous  with respect to $t$.

 (f2) $f(x, t) \leq C_0(1+t^p)$  uniformly for all $x$ in $\Omega$, and
 $$
 \mathop{\lim}\limits_{t\to \infty}\frac{f(x, t)}{t^p}=K(x), \,\, 1<p<\frac{n+2s}{n-2s},
 $$
 where $K(x)\in (0, \infty)$ is uniformly continuous and $ \mathop{\lim}\limits_{|x|\to \infty} K(x)=\bar C$.

\begin{rem}

(i) When $f(u)=u^p$, conditions (f1) and (f2) are satisfied.

\end{rem}

\begin{thm}\label{Ath1.5}
Let $\Omega$ be a bounded  domain in $\mathbb{R}^n$ and $u\in C_{loc}^{1, 1}(\Omega)\cap \mathcal{L}_{2s}$ is a nonnegative solution of \eqref{Aeq1.7}. Suppose that the assumptions (f1) and (f2) hold.

Then there exists a positive constant $C$ independent of solution $u$, such that
\be \label{Aest}
u(x) \leq C,\,\, x\in \Omega.
\ee
\end{thm}

\begin{rem}
As one can see from the following proof, the result of the theorem can be applied to more general equation \eqref{Aeq1.7}, such as
$$
\left\{\begin{array}{ll}
 \mathcal{L} u(x) + c(x) u(x) = f(x, u(x)), & x\in \Omega,\\
 u(x) \equiv 0, & x \in \Omega^c.
 \end{array}
 \right.
 $$

 In the case $1/2 <s < 1$, the left hand side of the equation can also contain $|\nabla u|^q$ for some proper choice of $q$.

 \end{rem}

We will use the blow-up and rescaling argument to obtain such a priori estimate.

{\em Outline of proof of Theorem \ref{Ath1.5}}

Suppose (\ref{Aest}) does not hold, then there exist a sequence
 of solutions $\{u_k\}$ to (\ref{Aeq1.7}) and a sequence of points $\{x^k\} \subset \Omega$ such that
\be
u_k(x^k)=\underset{\Omega}{\max}\,u_k:=m_k \ra \infty.
\ee

Let
\be
\lambda_k=m_k^{\frac{1-p}{2s}} \mbox{ and } v_k(x)=\frac{1}{m_k}u_k(\lambda_k x+x^k),\label{pe2.1}
\ee
then we have

\begin{eqnarray}\label{eq4-8}
(-\Delta)^s v_k(x)=
\lambda_k^\frac{2sp}{p-1} f(\lambda_kx+x^k, \lambda_k^{-\frac{2s}{p-1}}v_k(x))
:=F_k(x, v_k(x))
\end{eqnarray}
in
$$\Omega_k:=\{x \in \mathbb{R}^n \mid x=\frac{y-x^k}{\lambda_k},\,
y \in \Omega\}. $$

Let $d_k=dist(x^k, \partial\Omega)$. According to how fast that  $x^k$ is approaching the boundary $\partial \Omega$ as compared to $\la_k$, there are three possible cases. We will use contradiction arguments to
exhaust all three possibilities.
\vspace{12pt}

\emph{Case 1.} $\underset{k \ra \infty}{\lim}\frac{d_k}{\lambda_k} =\infty$.
\vspace{12pt}

In this case, $x^k$ remains in the interior or approaches boundary ``slowly'', and the blow-up occurs essentially in the interior of $\Omega$. Then after rescaling,
$\Omega$ is dilated into entire space, that is
$$\Omega_k \ra \mathbb{R}^n \mbox{ as } k \ra \infty.$$
Applying the well-known H\"{o}lder and Schauder estimates for the fractional equations (for example, see \cite{CLM}), there exists a positive function $v$ such that as
 $k \ra \infty$,
\be
v_k(x) \ra v(x) \mbox{ and hence }  F_k (x, v_k(x)) \to K(x^0)v^p(x).
\label{pe4a}
\ee
Furthermore,
\be
 (-\lap)^{s}v_k(x)\ra(-\lap)^{s}v(x)\label{pe4}.
\ee
This way, we arrive at a solution in $\R^n$:
\be
(-\lap)^{s}v(x)=K(x^0)v^p(x), \quad x \in  \mathbb{R}^n. \label{pe5}
\ee

This contradicts  the well-known Liouville theorem which asserts
that equation
 (\ref{pe5}) admits no positive solution.

\emph{Case 2.} $\underset{k \ra \infty}{\lim}\frac{d_k}{\lambda_k} = \;C\;>\;0$.
\vspace{12pt}

In this case, $\Omega_k$ is dilated into an upper half space,
 $$ \Omega_k \ra R^n_{C}:=\{x_n\geq -C\mid x \in \R^n \} \mbox{ as }k \ra \infty.$$
 Similar to \emph{Case 1}, we arrive at a limiting problem
\be \left\{\begin{array}{ll}
(-\lap)^{s}v(x)=K(x^0)v^p(x), & x \in  \mathbb{R}^n_{C} \\
v(x) = 0, & x \in (\R^n_C)^c
\end{array}
\right.
\label{pe7}
\ee
with a positive solution $v$. This again contradicts the existing Liouville theorem
 \cite{CFY}.

\emph{Case 3.} $\underset{k \ra \infty}{\lim}\frac{d_k}{\lambda_k} = \;0$. We rule out this possibility by showing that
$\{v_k(x)\}$ is uniformly Holder continuous on the boundary upon the construction of a super solution.
\bigskip

\section{Regularity of solutions}
\medskip

\subsection{Regularity lifting} \,
\medskip

We recall a few classical approaches to enhance the regularity of weak solutions. Here for each theorem, we only provide an outline of proof, while more details arguments can be found in book \cite{CL1}.

\subsubsection{Bootstrap} \index{bootstrap method}

Assume that $u \in H^1 (\Omega)$ is a weak solution of
\begin{equation}
- \lap u = u^p (x) \; , \;\; x \in \Omega.
\label{bo1}
\end{equation}
By Sobolev embedding, $u \in L^{\frac{2n}{n-2}} (\Omega)$. Applying the equation once,
we obtain $u \in W^{2, \frac{2n}{p(n-2)}} (\Omega)$. Via the embedding, $u \in L^q(\Omega)$.
If
the power $p$ is less than the critical number $\frac{n+2}{n-2}$, then $q > \frac{2n}{n-2}$. Continuing this way,
the integral power $q$ of $u$ can be enhanced repeatedly until $W^{2, q/p} (\Omega)$ is embedded into  $\in C^{\alpha} (\Omega)$. Finally
the ``Schauder Estimate'' will lift $u$ to $C^{2+\alpha} (\Omega)$
and hence be a classical solution. We call this the  {\rm bootstrap
method}.

The above argument does not work in the critical case when $p = \frac{n+2}{n-2}$. To
deal with this situation, we introduce  {\em regularity lifting
methods} in the next subsections.

\subsubsection{Regularity lifting by contracting operators}

This method has been used extensively in various forms in the authors'
previous works. The essence of the approach is well known in the
analysis community.  However, the version we present
here contains some new developments. Essentially, it is based on the
following Regularity
Lifting Theorem.

Let $V$ be a Hausdorff topological vector space. Suppose there are two extended
norms (i.e,. the norm of an element in $V$ might be infinity) defined
on $V$,
$$\|\cdot\|_X, \;\|\cdot\|_Y: V \ra [0,\infty] .$$
Let
$$X:=\{v\in V: \|v\|_X<\infty\} \; \mbox{ and } \; Y:=\{v\in V: \|v\|_Y<\infty\}.$$
Assume that spaces $X$ and $Y$ are complete under the corresponding norms and the convergence in $X$
or in $Y$ implies the convergence in $V$.

\begin{thm} (Regularity Lifting I)\index{Regularity Lifting Theorem I}.
Let $T$ be a contracting map from $X$ into itself and from $Y$ into itself.
Assume that $f \in X$ and that there exists a function
$g \in Z := X\cap Y$ such that $f=Tf+g$ in $X$.
Then $f$ also belongs to $Z$.
\label{thmRL}
\end{thm}

\begin{rem}  In practice, we usually choose $V$ to be the space of distributions, and $X$ and $Y$ to be suitable function spaces, for instance, $X=L^p(\Omega)$ and
  $Y = W^{1,q}(\Omega)$. Starting from a function $f$ belonging to the lower
regularity space $X$; if we can verify that $T$ is a contraction from
$X$ to itself and from $Y$ to itself, then we can lift $f$ to the higher regularity space $$Z = X\cap Y = L^p(\Omega) \cap W^{1,q}(\Omega).$$
\end{rem}

\subsubsection{Applications to PDEs}

Consider the weak
solutions of the equation involving the critical exponent:
\begin{equation} - \lap u = u^{\frac{n+2}{n-2}}
\label{ap1u}
\end{equation}
in a smooth bounded domain $\Omega$.

Writing $u^{\frac{n+2}{n-2}}= u^{\frac{4}{n-2}} u \equiv a(x) u$
we may study the regularity of the weak solution to the
following more general equation
\begin{equation}
- \lap u = a(x) u + b(x) . \label{ap2u}
\end{equation}
\begin{thm}
Assume that $a(x), b(x) \in L^{\frac{n}{2}} (\Omega)$. Let $u$ be
any weak solution of equation (\ref{ap2u}) in $H_0^1 (\Omega)$. Then
$u$  is in $L^p (\Omega)$ for any $1 \leq p < \infty$.
\label{thmapde}
\end{thm}

Applying the above theorem to equation (\ref{ap1u}), we
first conclude that $u$ is in $L^q (\Omega)$ for any $1<q<\infty$,
then through the equation, $u$ is in $W^{2,\gamma}
(\Omega)$ for any $1< \gamma < \infty$. . This implies that $u \in C^{1,\alpha} (\Omega)$ via
Sobolev embedding. Finally, by repeated applications of
Schauder estimates, we derive that $u \in C^{\infty} (\Omega)$.

\begin{cor}
If $u$ is a $H_0^1 (\Omega)$ weak solution of equation (\ref{ap1u}),
then $u \in C^{\infty} (\Omega)$ and hence is a classical solution.
\end{cor}

{\em Outline of proof of Theorem \ref{thmapde}}.

Let $G(x,y)$ be Green's function of $- \lap$ in $\Omega$. For a positive number $A$, define
$$ a_A (x) = \left\{ \begin{array}{ll}
a(x) & \mbox{ if } |a(x)| \geq A \\
0 & \mbox{ otherwise } ,
\end{array}
\right.
$$
and $$ a_B (x) = a(x) - a_A (x) .$$

Let $$ (T_A u)(x) = \int_{\Omega} G(x, y) a_A (y) u(y) d y .$$ Then
equation (\ref{ap2u}) can be written as
$$
u(x) = (T_A u)(x) + F_A (x) ,
$$
where
$$ F_A (x) = \int_{\Omega} G(x, y) [ a_B (y) u(y) + b(y) ] d y .$$

We  show that, for any $1 \leq p < \infty$,

i) $T_A$ is a contracting operator from $L^p (\Omega)$ to $L^p
(\Omega)$ for $A$ large, and

ii) $F_A (x)$ is in $L^p (\Omega)$.

Then by the ``Regularity Lifting Theorem'', we derive immediately
that $u \in L^p (\Omega)$.

\begin{rem}

From the proof, one can see that if $-\lap$ is replaced by the fractional Laplacian $(-\lap)^s$, then the conclusion of Theorem \ref{thmapde} is still valid.
\end{rem}

The {\em regularity lifting theorem} can also be applied to a general second order uniformly elliptic operator in divergence
form
\begin{equation}
L u = - \sum_{i, j=1}^n \left( a_{ij} (x) u_{x_i} \right)_{x_j} +
\sum_{i=1}^n b_i (x) u_{x_i} + c(x) u \label{Ldiv1}
\end{equation}
to lift the regularity of weak solutions from $W_0^{1, p} (\Omega)$ to $W^{2,p}(\Omega)$ under appropriate conditions.

\subsubsection{Regularity lifting by combinations of contracting and shrinking operators}
\,
\medskip

In the previous subsections, we introduced a Regularity Lifting Theorem and
its applications, in which, in order to lift the regularity of a solution
from a lower to a higher space (in terms of regularity), we required
that operator $T$ be contracting in \emph{both} spaces.
However, for a nonlinear operator in certain spaces, it is sometimes
very difficult, if not impossible, to prove it to be contracting---although
one may still be able to show that it is ``shrinking'' as we will defined in the theorem  below. Here we introduce
a more general theorem, which requires that the operator be
contracting in one space but only ``shrinking'' in the other.
We believe this theorem will find increasing applications
in many situations in nonlinear analysis.

Let $V$ be a Hausdorff topological vector
space. Suppose there are two extended norms (i.e., the norm of an element
in $V$ might be infinity) defined on $V$,
$$\|\cdot\|_X, \;\|\cdot\|_Y: V \ra [0,\infty] .$$
Let
$$X:=\{v\in V: \|v\|_X<\infty\} \; \mbox{ and } \; Y:=\{v\in V: \|v\|_Y<\infty\}.$$
We also assume that $X$ is complete and that the topology in
$V$ is weaker than both the topology of $X$ and the weak topology of $Y$. In other words, convergence in $X$, or weak convergence in $Y$, implies
convergence in $V$.

{\em The pair of spaces $(X,Y)$ described above is called an ``XY-pair'', if whenever the sequence
$\{u_n\}\subset X$ with $u_n\rightarrow u$ in $X$ and $\|u_n\|_Y\leq
C$ will imply $u\in Y$.}

In practice, we usually choose $V$ to be the space of distributions and $X,
Y$ to be suitable function spaces, such as $L^p$ spaces, H\"{o}lder spaces, Sobolev
spaces, and so on. There are many commonly used function spaces that
are ``XY-pairs,'' as will be illustrated by some examples in the remark following the theorem.

\begin{thm} \label{regularity-lifting}
(Regularity Lifting II.) Suppose Banach spaces $X,Y$ are an ``XY-pair'',
both contained  in some larger topological space $V$ satisfying
properties described above. Let $\mathfrak{X}$ and $\mathfrak{Y}$ be closed subsets of $X$ and $Y$ respectively. Suppose $T: \mathfrak{X} \rightarrow X$ is a contraction:
\begin{equation}
\|T f - T g\|_X \leq \eta \|f -g\|_X , \; \forall \, f, g \in \mathfrak{X} \; \text{ and for some } 0 < \eta <1 ;
\label{s5}
\end{equation}
and $T: \mathfrak{Y} \rightarrow Y$ is shrinking: \index{shrinking operator}
\begin{equation} \|T g\|_Y \leq \theta\|g\|_Y , \;  \forall \, g \in \mathfrak{Y}, \; \text{ and for some } 0<\theta<1.
\label{s6}
\end{equation}

Define
$$ S f = T f + F \;\;\; \mbox{ for some } F \in \mathfrak{X} \cap \mathfrak{Y} .$$
Moreover, assume that
\begin{equation} S: \mathfrak{X} \cap \mathfrak{Y} \ra \mathfrak{X} \cap \mathfrak{Y} .
\label{s7}
\end{equation}

Then there exists a unique solution $u$ of equation $$ u = T u + F  \; \mbox{ in } \; \mathfrak{X},$$ and more importantly,
$$ u \in Y .$$
\end{thm}

\begin{rem}  In some situations, one can
choose $\mathfrak{X} = X$ and $\mathfrak{Y} = Y$.

In practice, if one knows that a solution $u$ of $S u = u$ belongs to
$X$ (usually with lower regularity), then by Theorem \ref{regularity-lifting},
one can lift the regularity of $u$ up to $u \in X\cap Y$ (with higher
regularity).
\end{rem}

\begin{rem}
``XY-pairs'' are quite common, for examples
\begin{itemize}

\item $X=L^p(U)$ for $1\leq p\leq\infty$, $Y=C^{0,\alpha}(U)$ for
$0<\alpha\leq 1$, and $V$ is the space of distributions. Here $U$ can be
any subset of $\R^n$ or $\R^n$ itself.

\item $X$ is a Banach space, $Y$ is a reflexive Banach space, and
both are in some larger topological space $V$. Of course, we assume $V$ is
Hausdorff and has topology weaker than the topology of $X$ and the weak
topology of $Y$. Then for any $u_n\rightarrow u\in X$ and
$\|u_n\|_Y\leq C$, we have $u\in Y$.

Notice that all Hilbert spaces, such as $L^2$, $H^1$, and $H^2$,
are reflexive Banach spaces.
\end{itemize}
\label{rems1}
\end{rem}

\subsubsection{Applications to Integral Equations} \,
\medskip

Now, as a simple example to illustrate the idea, we apply the
Regularity Lifting Theorem II to the following integral equation
\begin{equation}
u(x) = \int_{\R^n} \frac{1}{|x-y|^{n-\alpha}} u^{\tau} (y) d y , \; x \in \R^n
\label{s1}
\end{equation}
with $\tau > \frac{n}{n-\alpha}.$  We will lift the regularity of any given positive solution from $L^{\infty}(\R^n)$ to
$C^{0,1}(\R^n)$--the space of Lipschitz continuous functions.

In the subsection preceding the previous, we used the Regularity Lifting
by Contracting Operators to show that
every positive solution $u$ of (\ref{s1}) is in $L^p (\R^n)$ for any $p > \frac{n}{n-\alpha}$. From here, one can derive immediately that
\begin{thm}
Every positive solution of equation (\ref{s1}) is uniformly bounded in $R^n$.
\label{thms2}
\end{thm}

\leftline{{\em Outline of proof.}}
\smallskip

Divide the integral into two parts
\begin{eqnarray*}
u(x) = \int_{B_1(x)} \frac{1}{|x-y|^{n-\alpha}} u^{\tau} (y) d y + \int_{\R^n \setminus B_1(x)} \frac{1}{|x-y|^{n-\alpha}} u^{\tau} (y) d y. \\
\end{eqnarray*}
Use the H\"{o}lder inequality on the first part is uniformly bounded, and it is obvious that the second part is uniformly bounded.

Based on Theorem \ref{thms2}, we will use Regularity Lifting Theorem II to lift the regularity of $u$ from $L^{\infty} (R^n)$ to $C^{0,1}(R^n)$, the space of Lipschitz continuous functions with norm
$$ \|v\|_{C^{0,1}(R^n)} = \|v\|_{L^{\infty} (R^n)} + \sup_{x \neq y} \left\{
\frac{|v(x) - v(y)|}{|x - y|} \right\} .$$

\begin{thm}
Every positive solution $u$ of equation (\ref{s1}) is Lipschitz continuous.
\label{thms3}
\end{thm}

\leftline{{\em Outline of proof}}
\smallskip

By elementary calculus, one can verify that
\begin{eqnarray*}
\int_0^{\infty} \int_{B_t(x)} u^{\tau} (y) d y \, \frac{d t}{t^{n-\alpha+1}} = (n-\alpha) \int_{R^n} \frac{1}{|x-y|^{n-\alpha}} u^{\tau} (y) d y .
\end{eqnarray*}
It follows that the solution of (\ref{s1}) only differs by a constant
multiple from the solution of the following equation:
\begin{equation}
u(x) = \int_0^{\infty} \int_{B_t(x)} u^{\tau} (y) d y \, \frac{d t}{t^{n-\alpha+1}} .
\label{s2}
\end{equation}
Hence, without loss of generality, we may just  prove that every positive solution
$u$ of (\ref{s2}) is Lipschitz continuous.

Let $$\mathfrak{X} = \{ v \in X \equiv L^{\infty} (R^n) \mid  \|v\|_{L^{\infty}} \leq 2 \|u\|_{L^{\infty}} \} $$
and $$\mathfrak{Y} = \{ v \in Y \equiv C^{0,1} (R^n) \mid \|v\|_{L^{\infty}} \leq 2 \|u\|_{L^{\infty}} \} .$$

For every $\epsilon > 0$, define
$$ T_{\epsilon} v (x) = \int_0^{\epsilon} \int_{B_t(x)} v^{\tau} (y) d y \, \frac{d t}{t^{n-\alpha+1}} $$
and $$ F(x) = \int_{\epsilon}^{\infty} \int_{B_t(x)} u^{\tau} (y) d y \, \frac{d t}{t^{n-\alpha+1}} .$$
Then obviously, $u$ is a solution of the equation
\begin{equation}
v = T_{\epsilon} v + F .
\label{s3}
\end{equation}
Write $S_{\epsilon} v = T_{\epsilon} v + F .$
We show that for $\epsilon$ sufficiently small,

(i) $T_{\epsilon}$ is a contracting operator from $\mathfrak{X}$ to $X$,

(ii) $T_{\epsilon}$ is a shrinking operator from $\mathfrak{Y}$ to $Y$,

(iii) $F \in \mathfrak{X} \cap \mathfrak{Y}$,  and $S_{\epsilon} : \mathfrak{X} \cap \mathfrak{Y} \ra \mathfrak{X} \cap \mathfrak{Y}$.
\medskip

\begin{rem}

From the proof (see \cite{CL1}), one can notice that it is extremely difficult, if not impossible,  to prove that $T_{\epsilon}$ is a contracting  operator from $\mathfrak{Y}$ to $Y$ due to
the complexity of the Lipschitz norm in $Y$.

\end{rem}

\subsubsection{Applications to fully nonlinear systems of Wolff type} \index{fully nonlinear Wolff systems} \,
\medskip

The Regularity Lifting Theorem II can also be applied to the following fully nonlinear integral systems involving Wolff potentials to prove that every positive solution pair is Lipschitz continuous (See \cite{CLM}). Consider
\begin{equation}
\left\{\begin{array}{ll}
u(x) = W_{\beta, \gamma}(v^q)(x) , & x \in \R^n,  \\
v(x) = W_{\beta, \gamma} (u^p)(x) , & x \in \R^n;
\end{array}
\right.
\label{WFPs}
\end{equation}
where
$$ W_{\beta,\gamma} (f)(x) = \int_0^{\infty} \left[ \frac{\int_{B_t(x)} f(y) dy}{t^{n-\beta\gamma}} \right]^{\frac{1}{\gamma-1}} \frac{d t}{t}.$$

This system includes many known systems as special cases, in particular, when $\beta = \frac{\alpha}{2}$ and $\gamma = 2$, system (\ref{WFPs}) reduces to
\begin{equation}
 \left \{
 \begin{array}{l}
      u(x) = \int_{\R^{n}} \frac{1}{|x-y|^{n-\alpha}} v(y)^q dy , \;\; x \in \R^n , \\
      v(x) = \int_{\R^{n}} \frac{1}{|x-y|^{n-\alpha}} u(y)^p dy , \;\; x \in \R^n .
\end{array}
\right. \label{systems1}
\end{equation}
The solutions $(u,v)$ of (\ref{systems1}) are critical points of the functional associated with
the well-known Hardy-Littlewood-Sobolev inequality. We can show that
(\ref{systems1}) is equivalent to a system of semi-linear elliptic PDEs
$$
 \left \{
 \begin{array}{l}
      (-\Delta)^{\alpha/2} u = v^q , \;\; \mbox{ in } \R^n,\\
      (-\Delta)^{\alpha/2} v = u^p, \;\;  \mbox{ in } \R^n .
\end{array}
\right.  $$
When $\alpha =2$, it becomes the well-known Lane-Emden system.

System (\ref{WFPs}) is said to be in critical case if
\begin{equation}
\frac{\gamma -1}{p+\gamma -1} + \frac{\gamma -1}{q +\gamma -1} = \frac{n-\beta \gamma}{n} .
\label{cri}
\end{equation}

In \cite{CLM}, the authors proved the following:

\begin{thm}
Let $(u,v)\in L^{p+\gamma-1}(\R^n)\times
L^{q+\gamma-1}(\R^n)$ be a pair of positive solutions for system (\ref{WFPs}) in the critical case (\ref{cri}).
Further assume $p>1$, $q>1$, and $1<\gamma \leq 2$. Then $(u,v)$ is uniformly bounded.
\label{pros1}
\end{thm}

Based on this theorem, the authors used Regularity Lifting Theorem II to
lift the positive solution pair from
$L^{\infty} (\R^n)$ to $C^{0,1}(\R^n)$.

\begin{thm}
Under the same conditions as in Theorem \ref{pros1}, $(u,v)$ is
Lipschitz continuous.
\label{thms5}
\end{thm}
\bigskip

\subsection{Interior regularity--elliptic case} \,
\medskip

Among other applications, interior regularity estimates are powerful tools in the blow-up and rescaling analysis to drive a priori estimates for solutions, as we will explain in a moment.

Consider
\be \label{RE}
(-\lap)^s u(x) = f(x), \;\; x \in B_1(0).
\ee

Currently available typical regularity estimates look like:

\begin{pro} (H\"{o}lder estimate)

Assume that $u$ is a solution of \eqref{RE} and $f \in L^\infty(B_1(0))$. Then for any $0<\varepsilon < 2s$, there exists a constant $C$, such that
\be \label{REH}
\|u\|_{C^{[2s-\varepsilon], \{2s-\varepsilon\}}(B_{1/2}(0))}\leq C(\|f\|_{L^\infty(B_1(0))}+\|u\|_{L^\infty(\mathbb{R}^n)}).
\ee
\end{pro}

\begin{pro} (Schauder estimate)
Assume that $u$ is a solution of \eqref{RE} and $f \in C^\alpha(\bar B_1(0))$. Then there exists a constant $C$, such that
\be \label{RES}
\|u\|_{C^{[2s+\alpha], \{2s+\alpha\}}(B_{1/2}(0))}\leq
C\left(\|f\|_{C^\alpha(\bar B_1(0))} +\|u\|_{L^\infty(\mathbb{R}^n)}\right),
\ee
where $2s+\alpha \notin \mathbb{N},$  $[2s+\alpha]$ and $\{2s+\alpha\}$ are  the integer part and the fraction part of $2s+\alpha$ respectively.
\end{pro}
These results can be found in  \cite{CLM, RS, Si}  and also in \cite{CS1, CS2}
for fully nonlinear nonlocal equations. In particular, in \cite{LW}, the authors established the pointwise regularity of the solutions for the fractional equation, from which they derived classical Schauder regularity.

Notice that, due to the nonlocality of the fractional Laplacian, the right hand side of both \eqref{REH} and \eqref{RES} contains $\|u\|_{L^\infty(\mathbb{R}^n)}$.  These known regularity estimates are insufficient  for applying the
powerful blow-up  and rescaling arguments to establish a priori estimates for fractional equations on {\em unbounded} domains due to the requirement on the bonded-ness
of global norms of $u.$

To illustrate such an insufficiency, let $\Omega$ be an {\em unbounded} domain in $\mathbb{R}^n$. We compare a simple problem of local nature
\begin{equation}
\left\{\begin{array}{ll}
- \Delta u(x) = u^p(x), & x \in \Omega,\\
u(x) = 0, & x \in \partial \Omega,
\end{array}
\right.
\label{1001}
\end{equation}
with its nonlocal counter part
\begin{equation}
\left \{\begin{array}{ll}
(- \Delta)^s u(x) = u^p(x), & x \in \Omega,\\
u(x) = 0, & x \in \mathbb{R}^n \setminus \Omega,
\end{array}
\right.
\label{1002}
\end{equation}
where $0<s<1$, and $p$ is a subcritical exponent.

To obtain an a priori estimate of positive solutions, one would usually apply a typical blowing up and re-scaling argument briefly as follows.

If the solutions are not uniformly bounded, then there exist a sequence of solutions $\{u_k\}$ and a sequence of points $\{x^k\} \in \Omega$,
such that $u_k(x^k) \to \infty$. Notice that since $\Omega$ is unbounded, this $u_k(x^k)$ may not be the maximum of $u_k$ in $\Omega$.
Upon re-scaling, the new sequence of functions $v_k$ are bounded and satisfy the same equation on certain domains $B_k$.

{\em For the  equation in \eqref{1001} of local nature}, the bounded-ness of $\{v_k\}$ on $B_k$  is  {\em sufficient} to guarantee the bound on its higher norms, say $C^{2, \alpha}$ norm.  This leads to the convergence of  $\{v_k\}$ to a function $v$, which is a bounded solution of the same equation as \eqref{1001} either in the whole space or in a half space. By the known results on the nonexistence of solutions in these two spaces, one derives a contradiction and hence obtains the a priori estimate.

{\em Nevertheless, when dealing with nonlocal problem  \eqref{1002} }, the currently available regularity estimates, such as   \eqref{REH} and \eqref{RES}, indicate that the bounded-ness of $\{v_k\}$ on $B_k$ are {\em not sufficient} to guarantee the bound on its higher norms, instead,  $\{v_k\}$ is required to {\em be bounded across the entire
space.} Such a condition cannot be met because there is no  information on the behavior of $v_k$ beyond $B_k$ during the re-scaling process. This posses  a substantial difficulty on employing blowing and re-scaling arguments for nonlocal equations on unbounded domains.
\medskip

The above arguments lead a natural and challenging question:

\medskip

Can one use only local bound $\|u\|_{L^\infty(B_1(0))}$ to replace the global bound $\|u\|_{L^\infty(\mathbb{R}^n)}$
 in \eqref{REH} and \eqref{RES}?
\medskip

The answer to this question is affirmative for nonnegative solutions as provided in \cite{CLWX}, in which the authors established
stronger interior H\"{o}lder and Schauder regularity estimates.

\begin{thm} \cite{CLWX} (H\"{o}lder estimate) \label{mthm6}

Assume that $u$ is a nonnegative solution of \eqref{RE} and $f \in L^\infty(B_1(0))$. Then for any $0<\varepsilon < 2s$, there exists a constant $C$, such that
\be \label{REH1}
\|u\|_{C^{[2s-\varepsilon], \{2s-\varepsilon\}}(B_{1/2}(0))}\leq C(\|f\|_{L^\infty(B_1(0))}+\|u\|_{L^\infty(B_1(0))}).
\ee
\end{thm}

\begin{thm} \cite{CLWX} (Schauder estimate) \label{mthm7}
Assume that $u$ is a nonnegative solution of \eqref{RE} and $f \in C^\alpha(\bar B_1(0))$. Then there exists a constant $C$, such that
$$
\|u\|_{C^{[2s+\alpha], \{2s+\alpha\}}(B_{1/2}(0))}\leq
C\left(\|f\|_{C^\alpha(\bar B_1(0))} +\|u\|_{L^\infty(B_1(0))}\right),
$$
where $2s+\alpha \notin \mathbb{N},$  $[2s+\alpha]$ and $\{2s+\alpha\}$ are  the integer part and the fraction part of $2s+\alpha$ respectively.
\end{thm}

The key improvement here is that the local bound $\|u\|_{L^\infty(B_1(0))}$ is employed instead of the global bound $\|u\|_{L^\infty(\R^n)}$ to
control the higher derivative norms. This advancement makes it possible to carry out the blow-up and rescaling analysis to establish a priori estimates for nonnegative solutions on unbounded domains as shown in the following example.
\medskip

Consider
\be \label{AP2}
\left\{\begin{array}{ll} (-\lap)^s u(x) + b(x) |\nabla u(x)|^q = f(x, u(x)), & x \in \Omega,\\
u(x) = 0, & x \in \Omega^c.
\end{array}
\right.
\ee

\begin{thm} \cite{CLWX} \label{mthm9} Suppose $s>\frac{1}{2}$
Let $\Omega$ be an arbitrary domain in $\mathbb{R}^n$ and $u\in C_{loc}^{1, 1}(\Omega)\cap \mathcal{L}_{2s}$ is a nonnegative solution of \eqref{AP2}. Suppose that the assumptions (f1), (f2) (see Section 4.4) hold and

 $b(x):  \mathbb{R}^n \to \mathbb{R}$  is uniformly bounded and uniformly H\"{o}lder continuous.

If $s>\frac{1}{2},$ then there exists a positive constant $C$ such that
$$
u(x)+|\nabla u|^{\frac{2s}{2s+p-1}}(x)\leq C\left(1+dist^{-\frac{2s}{p-1}}(x, \partial \Omega)\right),\,\, x\in \Omega.
$$
More precisely, if $\Omega \neq \mathbb{R}^n,$ we have
$$
u(x)+|\nabla u|^{\frac{2s}{2s+p-1}}(x)\leq C\, dist^{-\frac{2s}{p-1}}(x, \partial \Omega) ,\,\, x\in \Omega.
$$
In particular, if $\Omega$ is the whole space, i.e. $\Omega=\mathbb{R}^n,$ then
$$
u(x)+|\nabla u|^{\frac{2s}{2s+p-1}}(x)\leq C,\,\, x\in \mathbb{R}^n;
$$
if $\Omega \neq \mathbb{R}^n$ is an exterior domain, i.e. $\Omega \supset \{x\in \mathbb{R}^n \mid |x|>R\}$,   then
$$
u(x)+|\nabla u|^{\frac{2s}{2s+p-1}}(x)\leq C |x|^{-\frac{2s}{p-1}},\,\, |x|\geq 2R;
$$
and if $\Omega$ is a punctured ball, i.e. $\Omega=B_R(0)\backslash \{0\}$ for some $R>0,$ then
$$
u(x)+|\nabla u|^{\frac{2s}{2s+p-1}}(x)\leq C |x|^{-\frac{2s}{p-1}},\,\, 0<|x|\leq \frac{R}{2}.
$$
\end{thm}

To derive the Schauder estimates for non-negative solutions of the fractional Poisson equation \eqref{RE},
we split the solution $u$ into two parts: the potential
$$
w(x):=C_{n, s}\int_{B_1(0)}\frac{f(y)}{|x-y|^{n-2s}}dy.
$$
and the $s$-harmonic function
$$
h(x):=u(x)-w(x).
$$

The estimates for the potential $w(x)$ are rather standard.

For the $s$-harmonic function
$h(x)$, we express it in terms of the fractional Poisson kernel
\begin{equation*}
h(x)=\int_{|y|>1} P(x, y)h(y)dy, \,\, \forall \; |x|<1,
\end{equation*}
where
\begin{eqnarray*}
P(x, y)=
\begin{cases}
\frac{\Gamma (n/2)\sin(\pi s)}{\pi^{\frac{n}{2}+1}}\left(\frac{1-|x|^2}{|y|^2-1}\right)^s\frac{1}{|x-y|^n}, &|y|>1,\\
0,& |y|\leq 1.
\end{cases}
\end{eqnarray*}

By the definition of  $P(x, y)$ and some elaborate calculations, we find that
 the higher order derivatives of $P(x, y)$ can be controlled by itself. Then,
  $\|D^k h\|_{L^\infty(B_{1/2}(0))}$ can be bounded by $\|f\|_{L^\infty(B_1(0))}$ and  $\|u\|_{L^\infty(B_1(0))}$ due to the fact that $h=u-w$.
 This leads to better  Schauder estimates for nonnegative solutions.

\bigskip

\subsection{Interior regularity--parabolic case} \,
\medskip

Consider the fully fractional heat equation--the master equation
\begin{equation}\label{model}
    (\partial_t-\Delta)^s u(x,t)=f(x, t, u(x,t)) ,\,\,  \mbox{in}\,\,  \mathbb{R}^n\times\mathbb{R},
\end{equation}
where
\begin{equation}\label{nonlocaloper}
(\partial_t-\Delta)^s u(x,t)
:=C_{n,s}\int_{-\infty}^{t}\int_{\mathbb{R}^n}
  \frac{u(x,t)-u(y,\tau)}{(t-\tau)^{\frac{n}{2}+1+s}}e^{-\frac{|x-y|^2}{4(t-\tau)}}\operatorname{d}\!y\operatorname{d}\!\tau,
\end{equation}
 Note that this operator is nonlocal both in space and time, since the value of $(\partial_t-\Delta)^s u$ at a given point $(x,t)$ depends on the values of $u$ in the whole $\mathbb{R}^n$ and on  all the past time before $t$.
The singular integral in \eqref{nonlocaloper} is well defined provided
 $$u(x,t)\in C^{2s+\epsilon,s+\epsilon}_{x,\, t,\, {\rm loc}}(\mathbb{R}^n\times\mathbb{R}) \cap \mathcal{L}(\mathbb{R}^n\times\mathbb{R})$$
for some $\varepsilon\in (0,1)$, where
the slowly increasing function space $\mathcal{L}(\mathbb{R}^n\times\mathbb{R})$ is defined by
$$ \mathcal{L}(\mathbb{R}^n\times\mathbb{R}):=\left\{u \mid \int_{-\infty}^t \int_{\mathbb{R}^n} \frac{|u(x, \tau)|e^{-\frac{|x|^2}{4(t-\tau)}}}{1+(t-\tau)^{\frac{n}{2}+1+s}}\operatorname{d}\!x\operatorname{d}\!\tau<\infty,\,\, \forall \,t\in\mathbb{R}\right\}. $$

It can be verified that when the space-time nonlocal operator $(\partial_t-\Delta)^s$ is applied to a function that only depends on the spatial variable $x$, it becomes the well-known fractional Laplacian in the sense that, for each $x \in \R^n$,
 \begin{equation*}
   (\partial_t-\Delta)^s u(x)=(-\Delta)^s u(x).
 \end{equation*}
While in the case $u=u(t)$, then
 \begin{equation*}
   (\partial_t-\Delta)^s u(t)=\partial_t^s u(t),
 \end{equation*}
where $\partial_t^s$ is the Marchaud fractional derivative of order $s$. Moreover, as $s\rightarrow 1$, the fractional powers of heat operator $(\partial_t-\Delta)^s$
reduces to the local heat operator $\partial_t-\Delta$.

The space-time nonlocal equation represented by \eqref{model} arises in various physical and biological phenomena, such as anomalous diffusion, chaotic dynamics, biological invasions, and so on.

In this section, we investigate the interior regularity estimate for nonnegative solutions to such equations.

Let
$$Q_r (x^o, t_o) = \{(x,t) \in \R^n \times \R \mid |x-x^o| < r, \; |t-t_o| < r^2 \}$$
be the parabolic cylinder of size $r$ centered at $(x^o, t_o)$.

Regarding regularity estimates for master equation
\begin{equation}\label{ME1}
  (\partial_t-\Delta)^s u(x,t) =f(x,t),\,\,(x, t)\in \mathbb{R}^n\times \mathbb{R},
\end{equation}
Stinga and Torrea \cite{ST}, among other results, established the following local H\"{o}lder estimate ($0<s<1/2$ case)
\begin{equation}  \label{ST1}
\|u\|_{C_{x,t}^{2s, s}(Q_1(0,0))} \leq C \left( \|f\|_{L^\infty(Q_2(0,0))} + \|u\|_{L^\infty (\mathbf{\R^n \times (-\infty, 4)})}\right),
\end{equation}
and Schauder estimate:
\begin{equation}  \label{ST2}
\|u\|_{C_{x,t}^{2s + \alpha, s + \alpha/2}(Q_1(0,0))} \leq C \left( \|f\|_{C_{x,t}^{\alpha, \alpha/2}(Q_2(0,0))} + \|u\|_{L^\infty (\bf{\R^n \times (-\infty, 4)})}\right) .
\end{equation}

Here we use the standard notation for parabolic H\"{o}lder space. A function $u(x,t)$ is said to belong to $C_{x,t}^{\alpha, \beta}$ if $u$ is $C^\alpha$ in the spatial variable $x$ and $C^\beta$ in the time variable $t$. If $\alpha <1$, $C^\alpha$ is the usual H\"{o}lder space.
If $1<\alpha <2$, $C^\alpha$ denotes the space of bounded functions whose first derivatives belong to $C^{\alpha -1}$, and so on.

In order to circumvent the nonlocality of the master operator, the authors relied on the extension method introduced by Caffarelli and Silvestre \cite{CS}.
Due to the nonlocal nature of the fully fractional heat operator, the global norm $\|u\|_{L^\infty (\R^n \times (-\infty, 4))}$ needs to be placed on the right hand sides of \eqref{ST1} and \eqref{ST2}.

As we analyzed in the previous section, such an estimates are insufficient to apply the blow-up and rescaling argument to establish a priori estimate on the solutions in unbounded domains.

To answer this challenge, the authors in \cite{CGL1} replace the global norm $\|u\|_{L^\infty (\R^n \times (-\infty, 4))}$ with the local norm $\|u\|_{L^\infty (\mathbf{Q_2(0,0)}}$ when considering
nonnegative solutions $u(x,t)$ to master equation \eqref{ME1} with  $f(x,t)\geq 0$. The following are two typical estimates.

\begin{thm} \cite{CGL1} (H\"{o}lder estimates) \label{mthm14} Assume that both $f(x,t)$ and $u(x,t)$ are bounded in $Q_2(0,0)$. Then there exists a positive constant $C$, such that
$$
\|u\|_{C^{2s, s}_{x, t} (Q_1(0,0))} \leq C \, \left( \|f\|_{L^\infty (Q_2(0,0))} + \|u\|_{L^\infty (\mathbf{Q_2(0,0)})}\right), \;\; \mbox{ if } s \neq 1/2.
$$
\end{thm}

\begin{thm} \cite{CGL1} (Schauder estimates) \label{mthm15} Assume that $u\in L^\infty(Q_2(0,0)),  f  \in C^{\alpha, \alpha/2}_{x,t}  (Q_2(0,0))$ for some $0 < \alpha < 1$. Then there exists a positive constant $C$, such that
$$
\|u\|_{C^{2s+ \alpha, s+\alpha/2}_{x, t} (Q_1(0,0))} \leq C \, \left(\|f\|_{C^{\alpha, \alpha/2}_{x,t}  (Q_2(0,0))} + \|u\|_{L^\infty  (\mathbf{Q_2(0,0)})} \right),
$$
if $2s +\alpha < 1$ or if $1< 2s +\alpha < 2$.
\end{thm}
Other cases were also considered in the same paper.

We believe that this is the first instance where the local higher norms of solutions to a {\em nonlocal} equation can be effectively controlled by their local norms. This breakthrough makes it possible to employ the blow-up and re-scaling arguments to
establish a priori estimates for solutions to a family of similar nonlocal parabolic equations in unbounded domains.

Based on Theorems \ref{mthm14} and \ref{mthm15}, employing blow-up and rescaling analysis, the authors obtained a priori estimates for positive solutions to nonlinear master equations.

Consider
\begin{equation}
(\partial_t - \Delta)^s u(x,t) = f(x, u(x,t)), \;\; (x,t) \in \mathbb{R}^n \times \mathbb{R}.
\label{AA101}
\end{equation}

\begin{thm} \cite{CGL1} \label{mthm16}
Assume that $u$ is a nonnegative solution of (\ref{AA101}) and that the following conditions hold:

(f1) $f(x, \tau)$ is uniformly H\"{o}lder continuous with respect to $x$ and continuous  with respect to $\tau$.

 (f2) There exists a constant $C_0 > 0$ such that
 $$f(x, \tau) \leq C_0(1+\tau^p) \;\;  \mbox{   uniformly for all $x$, and }$$
 $$
 \mathop{\lim}\limits_{\tau \to \infty}\frac{f(x, \tau)}{\tau^p}=K(x), \,\, 1< p < \frac{n+2}{n+2-2s},
 $$
 where $K(x)\in (0, \infty)$ is uniformly continuous and
 $$ \mathop{\lim}\limits_{|x|\to \infty} K(x)=\bar C\in (0, \infty).$$
 Then there exists a constant $C$, such that  for all nonnegative solutions $u$, we have
 $$
 u(x,t) \leq C, \,\,\, \forall \, (x,t) \in \mathbb{R}^n \times \mathbb{R}.
 $$

 \end{thm}
 The above a priori estimate holds for all $0<s<1$.

 In the case $1/2 <s <1$, they  further established estimates for more general master equations of the form
 \begin{equation} \label{102}
 (\partial_t - \Delta)^s u(x,t) = b(x)|\nabla_x u(x,t)|^q + f(x, u(x,t)), \;\; (x,t) \in \mathbb{R}^n \times \mathbb{R},
 \end{equation}
where $\nabla_x u$ is the gradient of $u$ with respect to $x$.

\begin{thm} \cite{CGL1} \label{mthm17} Assume that
  $b(x):  \mathbb{R}^n \to \mathbb{R}$  is uniformly bounded and uniformly H\"{o}lder continuous and $f(x, \tau)$ satisfies (f1) and (f2).

If $s>\frac{1}{2}$ and $0<q<\frac{2sp}{2s+p-1},$ then there exists a positive constant $C$ such that
$$
u(x,t)+|\nabla_x u|^{\frac{2s}{2s+p-1}}(x,t)\leq C,  \;\;\; \forall \, (x,t) \in \mathbb{R}^n \times \mathbb{R},
$$
for all nonnegative solutions $u(x,t)$ of \eqref{102}.
\end{thm}
\medskip

To established the regularity, the authors first prove that the pseudo-differential equation is equivalent to the integral equation
$$u(x,t) = c + \int_{-\infty}^t \int_{\mathbb{R}^n} f(y,\tau) G(x-y, t-\tau) dy d\tau,$$
where
$$G(x-y, t-\tau) = \frac{C_{n,s}}{(t-\tau)^{n/2+1-s}} e^{-\frac{|x-y|^2}{4(t-\tau)}}.$$

Since the constant $c$ does not influence either the H\"{o}lder norm or the derivatives of $u$, without loss of generality, we may consider
\begin{equation} \label{inteq1}
u(x,t) = \int_{-\infty}^t \int_{\mathbb{R}^n} f(y,\tau) G(x-y, t-\tau) dy d\tau,
\end{equation}

For simplicity of notation, we denote
$$ Q = Q_2(0,0), \mbox{ and } \; \; \tilde{Q} = Q_1(0,0) \subset \subset Q.$$

We split $u(x,t)$ into two parts $u(x,t)=v(x,t) + w(x,t)$ as  in the following.

Write
$$ f_Q (x,t) = \left\{ \begin{array}{ll} f(x,t), & (x,t) \in Q ,\\
0, & (x, t) \in Q^c.
\end{array}
\right.
$$

Denote
$$w(x,t) = \int_{-\infty}^t \int_{\mathbb{R}^n} f_Q(y,\tau) G(x-y, t-\tau) dy d\tau$$
Then $w(x,t)$ satisfies the nonhomogeneous equation in $Q$:
$$(\partial_t -\Delta)^s w(x,t) = f(x,t), \;\; \forall \, (x,t) \in Q.$$
While the homogeneous part can be expressed by
$$v(x,t):=u(x,t)-w(x,t)=\int_{-\infty}^t \int_{\mathbb{R}^n} f_{Q^c}(y,\tau) G(x-y, t-\tau) dy d\tau$$
with $f_{Q^c}:=f-f_Q$. Apparently, $v$ is a solution of the homogeneous equation in $Q$:
$$(\partial_t - \Delta)^s v(x,t) = 0, \;\; \forall \, (x,t) \in Q.$$
\smallskip

We estimate $v(x,t)$ and $w(x,t)$ separately. We first obtain
\begin{thm} \label{mthmv} Assume that $v$ is bounded in $Q$, then for any $\alpha \in (0,1)$,
$$
\|v\|_{C_{x,t}^{2, \alpha} (\tilde{Q})} \leq C \|v\|_{L^\infty (Q)}.
$$
\end{thm}

\begin{rem}
Using similar arguments as in the proofs (see Section 3 in \cite{CGL}), it can be shown that the $C_{x,t}^{2, \alpha}$ norm on the left hand side can be replaced by $C_{x,t}^{k, \alpha}$ norm for any integer $k$. Nonetheless, this $C_{x,t}^{2, \alpha}$ norm is sufficient for the  blow-up and rescaling analysis.
\end{rem}

\leftline{\em New ideas involved in the proof}
\smallskip

To estimate the derivatives of $v(x,t)$, say $v_{x_i}$,  in terms of $v(x,t)$ itself, we encounter a significant challenge.
From the definition
$$v(x,t) = C_{n,s}\int_{-\infty}^t \int_{\mathbb{R}^n}  \frac{f_{Q^c}(y,\tau)}{(t-\tau)^{n/2+1-s}} e^{-\frac{|x-y|^2}{4(t-\tau)}} dy d\tau, $$
it follows that
$$v_{x_i}(x,t)=C_{n,s}\int_{-\infty}^t \int_{\mathbb{R}^n}  \frac{f_{Q^c}(y,\tau)}{(t-\tau)^{n/2+1-s}} \frac{x_i-y_i}{2(t-\tau)}e^{-\frac{|x-y|^2}{4(t-\tau)}} dy d\tau.$$

Comparing these two expressions, it appears that, in order to control $v_{x_i}$ in terms of $v$, we need to establish the inequality
\begin{equation}\label{est-fuc1}\frac{x_i-y_i}{2(t-\tau)}e^{-\frac{|x-y|^2}{4(t-\tau)}}\leq Ce^{-\frac{|x-y|^2}{4(t-\tau)}}.\end{equation}
However, this seems impossible because the extra term $\frac{x_i-y_i}{2(t-\tau)}$ on the left hand side is clearly unbounded. To overcome this difficulty, we introduce an innovative idea. Instead of attempting to control $|v_{x_i}(x,t)|$ directly in terms of $v(x,t)$ at the same point $(x,t)$, we find finitely many nearby points $(x^j,t)$ to accomplish this task:
$$ |v_{x_i}(x,t)| \leq C \sum_j^{2^n} v(x^j, t).$$

To illustrate this approach, let us consider the case of two dimensions ($n=2$).

For each fixed  $x\in B_1(0)$, we take $x$ as the center and divide the plane $\R^2$ into four quadrants $I_1, I_2, I_3, I_4.$
On each quadrant $I_j$, choose the point $x^j$ as the intersection of $\partial B_{1/\sqrt{2}}(x)$ and the diagonal of $I_j, \,  j=1, \cdots, 4$:
$$ x^j = x + \eta_j, \mbox{ with } \eta_1 = \frac{1}{\sqrt{2}}(\cos \frac{\pi}{4}, \sin \frac{\pi}{4}), \eta_2 = \frac{1}{\sqrt{2}} (\cos \frac{3\pi}{4}, \sin \frac{3\pi}{4}), \cdots.$$
We verify that for all $y\in B_2^c(0)\cap I_j$ it holds
\be \label{AA2} |y-x|^2 \geq |y-x^j|^2 + \frac{1}{2} |y-x|. \ee

 Using this inequality, we are able to control the left hand side of \eqref{est-fuc1} at point $x$ by the right hand side at point $x^j$ {\bf for all} $ y \in B_2^c (0) \cap I_j$:
 $$
\frac{|y-x|}{t-\tau} e^{-\frac{|y-x|^2}{4(t-\tau)}} \leq C e^{-\frac{|y-x^j|^2}{4(t-\tau)}}.
$$

Consequently
for all $(x, t) \in \tilde{Q}$, we obtain
$$|v_{x_i} (x, t)| \leq C\sum_j^4  v(x^j,t) \leq C \|v\|_{L^\infty (Q)}. $$

The following Figure 4 shows the case of $n=2$. Here $\R^2$ is divided into four quadrants centered at $x$.  For $y$ in the third quadrant $I_3$, we choose $x^3$ as the intersection of $\partial B_{1/\sqrt{2}}(x)$ and the diagonal $l_3$. By applying the {\em law of consine} to the triangle $\triangle yx^3x$, one can easily verify \eqref{AA2} (when $j=3$).

\begin{center}
\begin{tikzpicture}
    \def\circleA{2} 
    \def\circleO{2}
    \def\circleB{4} 

    \draw[thick,-, blue] (-5,0) -- (6,0) node[right] {};
    \draw[thick,-,blue] (0,-5) -- (0,5) node[above] {};
    \draw[thick, blue] (0,0) circle (\circleA) node[xshift=-38pt, yshift=12pt] {$B_1(x)$};
    \draw[thick](1,-0.5) circle (\circleO)node[xshift=40pt, yshift=-10pt] {$B_1(o)$};
    \draw[thick] (1, -0.5) circle (\circleB)node[xshift=120pt, yshift=-50pt] {$B_2(o)$} ;

    \filldraw[black] (1, -0.5) circle (2pt) node[below] {$o$};
      \coordinate (x) at (0,0);
      \def\radius{2};
      \path[name path=circle] (x) circle (\radius);
      \path[name path=line1] (-4,-5) -- (4.2,5);
      \path[name path=line2] (4,-5) -- (-4,5);
      \path[name intersections={of=circle and line1, by={P1,P2}}];
      \path[name intersections={of=circle and line2, by={Q1,Q2}}];
    \filldraw[red] (P1) circle (2pt) node[xshift=-4pt, yshift=10pt] {${x}^1$};
    \filldraw[red] (P2) circle (2pt) node[xshift=4pt, yshift=-10pt] {${x}^3$};
    \filldraw[red] (Q1) circle (2pt) node[above] {${x}^2$};
    \filldraw[red] (Q2) circle (2pt) node[xshift=12pt] {${x}^4$};
    \filldraw[red] (-5,-2)  circle (2pt) node[below] {$y$};

    \filldraw (0,0)  circle (2pt) node[right] {$x$};

    \draw[thick, purple] (-3.3,-4) -- (3.5,4) node[above] {};
    \draw[thick, purple] (3.2,-4) -- (-3.2,4) node[above] {};

    \draw[thick, green] (-5,-2) -- (0,0) node[above] {} ;
    \draw[thick, green] (-5,-2) -- (P2) node[below] {} ;

    \node[blue] at (2.9,-4.5) {$I_4$};
    \node[blue] at (2.9, 4.5) {$I_1$};
    \node[blue] at (-2.9,4.5) {$I_2$};
    \node[blue] at (-2.9,-4.5) {$I_3$};
     \node[purple] at (3.5,-4) {$l_4$};
    \node[purple] at (3.5, 3.5) {$l_1$};
    \node[purple] at (-3.2,4.2) {$l_2$};
    \node[purple] at (-3.4,-4.2) {$l_3$};
     \node[purple] at (-0.23,-0.23) {$\theta$};
\node [below=1cm, align=flush center,text width=8cm] at (0,-4.1)
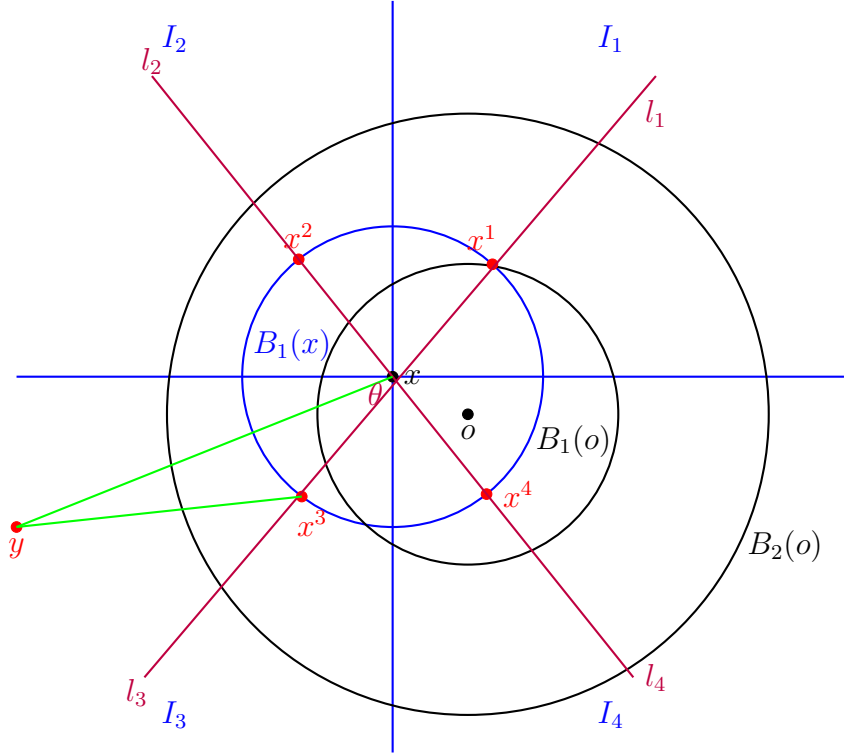
        {Figure 4.  {The case of $n=2$}.};

\end{tikzpicture}
\end{center}

The above approach can be interpreted as a {\em directional perturbation average}. At a given point $x$, it is impossible to control the derivative of $v$ solely by the value of $v$ itself. However, after making directional perturbations from $x$ to $x^j = x + \eta_j$ for $ j=1,2, \cdots, 2^n$ along directions $\eta_j$, and then taking the average, we are able to realize such a control. This estimate is essentially conducted on the fully fractional heat kernel $G(x-y, t-\tau)$. For convenience of future applications, we summarized it as
follows:
\begin{lem}[ Fully  fractional heat kernel estimates ]\label{key0} In each quadrant $I_j$ of $R^n$, there is a vector $\eta_j\in\partial B_{\frac{1}{\sqrt{n}}}(0)$ along the diagonal direction such that
\begin{equation}\label{key1} G(y,\tau)\leq {\bf{e^{-\frac{|y|}{n\tau}}}}\sum\limits_{j=1}^{2^n}G(y+{\eta}_j,\tau),  \forall \tau>0, y\in B_1^c(0).\end{equation}
Furthermore,   for each $\tau>0, y\in B_1^c(0)$, it holds
\begin{equation} D_y^\alpha  G(y,\tau)\leq C\sum\limits_{j=1}^{2^n}G(y+{\eta}_j,\tau), \forall \, |\alpha|=1 \, or \, 2,\end{equation}
and
 \begin{equation}  \partial_\tau G(y,\tau)\leq C\sum\limits_{j=1}^{2^n}G(y+{\eta}_j,\tau),\end{equation}
here $C=C(n)$ is a universal constant.
\end{lem}
 We believe that this important heat kernel estimate will become a powerful tool in establishing regularity estimates, with potential applications to a variety of related problems.
\medskip

\medskip

\subsection{Boundary regularity} \,
\medskip

Now consider the fractional Dirichlet problem in unbounded domains with boundary
\begin{equation} \label{A200}
        \begin{cases}
            \flap u(x) = g(x),&\text{in } \Omega,\\
            u(x)=0, &\text{in }  \Omega^c.\\
        \end{cases}
    \end{equation}

  As mentioned in Section 5.2, the authors in \cite{CLWX} established interior H\"older and Schauder regularity estimates for fractional equation \eqref{A200}. This makes it possible to derive a priori estimate for nonnegative solutions when $\Omega=\R^n$, an  unbounded domain {\em without boundary}.

Then for general unbounded domains  with  boundaries, can one establish a priori estimates for the solutions?

As we have seen in Section 4.4 when using the blow-up and rescaling argument to derive the a priori estimate, the sequence may blow-up on the boundary. To rule out this possibility, one needs to establish some boundary regularity, such as boundary H\"{o}lder continuity.

This requires boundary H\"older estimates in terms of only the local $L^{\infty}$ norm of the solutions.

However, in the existing results, such as the ones in \cite{BGQ} and \cite{BPGQ1}, in order to obtain the $C^s$ boundary H\"older continuity of the solutions for \eqref{A200} when $\Omega=\R^n_+$, the authors  require that solution $u$ is globally bounded.

With these existing boundary  regularity estimates, it is inadequate to carry out blow up and rescaling arguments on unbounded domains aimed at obtaining a priori estimates, because the rescaled functions may be globally unbounded.

In \cite{GLO}, the authors established the local version of the boundary regularity for nonnegative solutions, {\em only the local $L^{\infty}$  norm of $u$ instead of the global ones} is used to control its local H\"older norm up to the boundary.
\begin{thm}\label{1bdry C^s}
    Suppose $\Omega$ is a locally $C^{1,1}$ domain, $0<s<1$, $g\in L^\infty(\Omega\cap B_4)$ and $u$ is a nonnegative classical solution of
    \begin{equation}
        \begin{cases}
            \flap u = g(x),&\text{in } \Omega\cap B_4,\\
            u=0, &\text{in } B_4\backslash \Omega.\\
        \end{cases}
    \end{equation}
    Then $u\in C^{s}(\Omega\cap B_{1/2})$ and
    \[
    \norm{u}_{C^s(\Omega\cap B_{1/2})}\le C(\norm{g} _{L^\infty(\Omega\cap B_4)} + \norm{u} _{L^\infty(\Omega\cap B_4)} ),
    \]
    where $C$ depends on $n,s$ and $C^{1,1}_{loc}$ norm of $\partial\Omega$, and $B_r$ denotes a ball centered at any fixed point on the boundary $\partial\Omega$ with radius $r$.
\end{thm}

As an important application, they established a priori estimate for solutions  to
\begin{equation} \label{A201}
        \begin{cases}
            \flap u(x) = f(x,u(x)),&\text{in } \Omega,\\
            u(x)=0, &\text{in }  \Omega^c.\\
        \end{cases}
    \end{equation}

Assume that $\Omega \subset \R ^n$ is an unbounded domain with uniformly $C^{1,1}$ boundary, and $f$ satisfies the following condition:
\begin{itemize}
    \item $f(x,t):\Omega \times [0,\infty)\to\R$ is uniformly H\"older continuous with respect to $x$ and continuous with respect to $t$.
\end{itemize}
\begin{thm}\label{1A1}
Assume $1<p<\frac{n+2s}{n-2s}$, $f$ satisfies
\begin{equation}\label{growth f}
    f(x,t)\le C_0(1+t^p)
\end{equation}
uniformly for all $x$ in $\Omega$, and
\begin{equation}\label{lim f}
    \lim _{t\to\infty}\frac{f(x,t)}{t^p} = K(x),
\end{equation}
    where $K(x)\in (0,\infty)$ is uniformly continuous and $K(\infty):=\lim _{|x|\to\infty }K(x) \in (0,\infty)$.
Then there exists a constant $C$, such that
\begin{equation}\label{bd1}
u(x)\le C,\quad\forall \, x\in \Omega
\end{equation}
holds for all nonnegative solutions $u$ of \eqref{A201}.
\end{thm}
\bigskip
\bigskip

{\bf{Declaration:}}

Chen is
 partially supported by MPS Simons Foundation 847690  and  National Natural Science
Foundation of China (NSFC Grant No. 12571225)

 Guo is partially supported by  National Natural Science
Foundation of China (NSFC Grant No. 12501145), the Postdoctoral Fellowship Program of CPSF (No.GZC20252004),
the Natural Science Foundation of Shanghai (No.25ZR1402207), and the China Postdoctoral
Science Foundation (No.2025T180838 and 2025M773061).

Li is partially supported by NSFC-W2531006, NSFC-12250710674, NSFC-12031012 and the Institute of Modern Analysis-A Frontier Research Center of Shanghai.
\medskip

{\bf{Date availability statement:}} No data was used for the research described in the article.
\medskip

{\bf{Conflict of interest statement:}} There is no conflict of interest.

Wenxiong Chen

Department of Mathematical Sciences

Yeshiva University

New York, NY, 10033, USA

wchen@yu.edu
\medskip

Yahong Guo

School of Mathematical Sciences

Shanghai Jiao Tong University

Shanghai, 200240, P.R. China

yhguo@sjtu.edu.cn
\medskip

Congming Li

School of Mathematical Sciences

Shanghai Jiao Tong University

Shanghai, 200240, P.R. China

congming.li@sjtu.edu.cn


\begin{thebibliography}{CL}

\bibitem[ACV]{ACV} M. Allen, L. Caffarelli, and A. Vasseur,\,\, A parabolic problem with a fractional time derivative,\,\, Arch. Rational
    Mech. Anal., 221 (2016), 603-630.

\bibitem[AL]{AL} H. Amann and J. Lopez-Gomez, A priori bounds and multiple solutions for
superlinear indefinite elliptic problems, J. Diff. Eqns. 146 (1998),
336-374.

\bibitem[AT]{AT} S. Alama and G. Tarantello, On semilinear elliptic equations with indefinite nonlinearities,
Calc. Var. 1 (1993), 439-475.
\bibitem[B]{B} A.V. Balakrishnan, Fractional powers of closed operators and the semigroups generated by them, Pac. J. Math. 10 (1960) 419-437.

\bibitem[BCN]{BCN} H. Berestycki, I. Capuzzo-Dolcetta, and L. Nirenberg, Supperlinear
indefinite elliptic problems and nonlinear Liouville theorems,
Topol. Methods Nonl. Anal. 4 (1994), 59-78.

\bibitem[BCN1]{BCN1} H. Berestycki, I. Capuzzo-Dolcetta, and L. Nirenberg, Variational
methods for indifinite superlinear homogeneuous elliptic problems,
Nonlinear Diff. Eqns. Appl. 2 (1995), 553-572.



\bibitem[BBG]{BBG}  M. T. Barlow, R. F. Bass, and C. Gui, The Liouville property and a conjecture of De Giorgi, Comm. Pure Appl. Math., 53 (2000), 1007-1038.

\bibitem[BCN2]{BCN2} H. Berestycki, L. A. Caffarelli and L. Nirenberg, \,\, Monotonicity for elliptic  equations in unbounded Lipschitz domains, \,\, Comm. Pure Appl.  Math. 50(1997), 1089-1111.

\bibitem[BCPS] {BCPS} C. Br\"{a}ndle, E. Colorado, A. de Pablo, and U. S\'{a}nchez, \,\, A concave-convex elliptic problem involving the fractional Laplacian, \,\, Proc Royal Soc. of Edinburgh, 143 (2013), 39-71.

\bibitem[BGQ]{BGQ}  B. Barrios, J. Garcia-Melian, and A. Quaas. A note on the monotonicity of solutions for
fractional equations in half-spaces. Proc. Amer. Math. Soc., 147(7), (2019),3011-3019.



\bibitem[BHM]{BHM} H. Berestycki, F. Hamel and R. Monneau,\,\, One-dimensional symmetry of bounded entire solutions of some elliptic equations, \,\, Duke Math. J. 103(2000), 375-396.

\bibitem[BKN]{BKN} K. Bogdan, T. Kulczycki and A. Nowak, \,\, Gradient estimates for harmonic and $q$-harmonic functions of symmetric stable processes, \,\, Illinois J. Math.  46(2002) 541-556.




\bibitem[BN3]{BN3}  H. Berestycki and L. Nirenberg, \,\, On the method of moving planes and the sliding method, \,\, Bol. Soc. Brasil. Mat. (N.S.) 22(1991), 1-37.

\bibitem[BPGQ]{BPGQ} B. Barrios, L. Del Pezzo, J. Garcia-Melian and A. Quaas, \,\, A priori bounds and existence of solutions for some nonlocal elliptic problems, \,\, Revista Matematica Iberoamericana, 34(2018) 195-220.

\bibitem[BPGQ1]{BPGQ1} B. Barrios, L. Del Pezzo, J. Garcia-Melian, and A. Quaas, Symmetry results in the
half-space for a semi-linear fractional Laplace equation, Ann. Mat. Pura Appl. (4), 197(5):1385-1416, 2018.

\bibitem[BPSV]{BPSV} B. Barrios, I. Peral, F. Soria, and E. Valdinoci,\,\, A Widder's type theorem for the heat equation with nonlocal diffusion,  \,\,  Arch.   Rational     Mech. Anal.,  213 (2014),  629-650.

\bibitem[BV]{BV} C. Bucur, E. Valdinoci, \,\, Nonlocal diffusion and applications, \,\, Switzerland: Springer International Publishing, 2016.


\bibitem[CDL]{CDL}  W. Chen, L. D'Ambrosio, and Y. Li, \,\, Some Liouville theorems for the fractional Laplacian, \,\, Nonlinear Analysis, Theory, Methods \& Appl, 121(2015) 370-381.
\bibitem[CFS]{CFS} M. Caselli, E. Florit-Simon, J.  Serra,  Fractional Sobolev spaces on Riemannian manifolds. Math. Ann. 390, 6249-6314 (2024). https://doi.org/10.1007/s00208-024-02894-w
\bibitem[CFY]{CFY} W. Chen, Y. Fang, and R. Yang, \,\, Liouville theorems involving the fractional Laplacian on a half space, \,\, Adv. Math., 274 (2015), 167-198.

\bibitem[CG]{CG} W. Chen and Y. Guo, \,\, Master equations with indefinite nonlinearities, \,\, preprint 2024.

\bibitem[CGL]{CGL} W. Chen, Y. Guo, and C. Li, \,\, Lecture notes on the Fractional elliptic and parabolic equations, \, preprint, 2025.

\bibitem[CGL1]{CGL1} W. Chen, Y. Guo, and C. Li, \, Regularity of solutions for fully fractional parabolic equations, \,  2025. http://arxiv.org/abs/2502.07530

\bibitem[CGW]{CGW} W. Chen, Y. Guo, and L. Wu, Monotonicity for fractional semi-linear problem in a half space, preprint, 2025.

\bibitem[CH]{CH} W. Chen, Y. Hu,\,\, Monotonicity of positive solutions for nonlocal problems in unbounded domains,\,\, J. Funct. Anal., 281 (2021), 109187.


\bibitem[CHM]{CHM} W. Chen, Y. Hu, and L. Ma, \,\, Moving planes and sliding methods for fractional elliptic and parabolic equations, \, Adv. Nonlinear Stud. 24(2024), 359-398.


\bibitem[CL]{CL} W. Chen and C. Li, \,\, Classification of solutions of some nonlinear elliptic equations, Duke Math. J.  63(1991) 615-622.

\bibitem[CL1]{CL1} W. Chen, C. Li, \,\, Methods on nonlinear elliptic equations,\,\, AIMS book series, vol. 4, 2010.

\bibitem [CL2] {CL2} W. Chen, C. Li, \,\, Maximum principles for the fractional $p$-Laplacian and symmetry of solutions, \,\, Adv. Math., 335 (2018), 735-758.

\bibitem[CL3]{CL3} W. Chen and C. Li, \,\, A note on Kazdan-Warner type
conditions, \,\,  J. Diff. Geom.  41(1995), 259-268.

\bibitem[CL4]{CL4} W. Chen and C. Li, \,\,  A necessary and sufficient condition
for the Nirenberg problem,  \,\,   Comm.
 Pure Appl.  Math. 48(1995), 657-667.

\bibitem [CLg]{CLg} W. Chen, C. Li, and G. Li, \,\, Maximum principles for a fully nonlinear fractional order equation and symmetry of solutions, \,\, Calc. Var., 56 (2017), 29.


\bibitem[CLL]{CLL} W. Chen, C. Li, and Y. Li, \,\, A direct method of moving planes for the fractional Laplacian, \,\, Adv. Math., 308 (2017), 404-437.

\bibitem[CLL1]{CLL1}  W. Chen, C. Li, and Y. Li, \,\, A direct blow-up and rescaling argument on nonlocal elliptic equations, \,\,  International J. Math., 27 (2016), 1-20.

\bibitem [CLM]{CLM}  W. Chen, Y. Li, and P. Ma,\,\,  The fractional Laplacian,\,\, World Scientific, Hackensack, 2020.

\bibitem[CLO]{CLO} W. Chen, C. Li, and B. Ou, \,\, Classification of solutions for an integral equation, Commun. Pur. Appl. Math., 59 (2006), 330-343.

\bibitem[CLO1]{CLO1} W. Chen, C. Li, and B. Ou,  \,\,  Qualitative properties of solutions for an integral equation, \,\, Disc. Cont. Dyn. Sys., 12 (2005), 347-354.

\bibitem[CLO2]{CLO2} W. Chen, C. Li, and B. Ou, \,\, Classification of solutions for a system of integral equations, \,\, Commun. Partial Differ. Equ., 30 (2005), 59-65.

\bibitem[CLWX]{CLWX} W. Chen, C. Li, L. Wu, and Z. Xin, \, Refined regularity for nonlocal elliptic equations and applications,  2024. http://arxiv.org/abs/2502.04676


\bibitem[CLZ]{CLZ}  W. Chen, Y. Li, and R. Zhang, \,\, A direct method of moving spheres on fractional order equations, \,\, J. Funct. Anal., 272 (2017), 4133-4157.

\bibitem[CLZ1]{CLZ1} W. Chen, C. Li, and J. Zhu, \,\, Fractional equations with indefinite nonlinearities, \,\, Disc. Cont. Dyn. Sys., 39 (2019), 1257-1268.

\bibitem[CM]{CM} W. Chen, L. Ma, \,\, Qualitative properties of solutions for dual fractional nonlinear parabolic equations\,\, J. Funct. Anal. 285(2023),110117.

\bibitem[CM1]{CM1} W. Chen, L. Ma, \,\, Gibbon's conjecture for entire solutions of master equations, \, to appear in Comm. Compt. Math, 2024.

\bibitem[CS]{CS} L. Caffarelli, L. Silvestre, \,\, An extension problem related to the fractional Laplacian, \,\, Commun. Partial Differ. Equ., 32 (2007), 1245-1260.


\bibitem[CS1]{CS1} L. Caffarelli, L. Silvestre, \,\, A nonlocal Monge-Amp\`{e}re equation, \,\, Commun. Anal. Geom., 24 (2016), 307-335.

\bibitem[CS2]{CS2}
L. Caffarelli, L. Silvestre, The Evans-Krylov theorem for nonlocal fully nonlinear equations, Ann. of Math., 174 (2) (2011), 1163-1187.



\bibitem[CSt]{CSt} L. Caffarelli, P. Stinga, \,\, Fractional elliptic equations, Caccioppoli estimates and regularity, \,\, Ann. I. H. Poincar\'{e} - AN, 33 (2016), 767-807.

\bibitem[CV]{CV} L. Caffarelli, A. Vasseur, \,\, Drift diffusion equations with fractional diffusion and the quasi-geostrophic equation, \,\, Annals of Math., 171 (2010), 1903-1930.

\bibitem[CWNH]{CWNH} W. Chen, P. Wang, Y. Niu, and Y. Hu, \,\, Asymptotic method of moving planes for fractional  parabolic  equations, \,\, Adv. Math., 377 (2021), 107463.

\bibitem[CW]{CW} W. Chen, L. Wu, \,\, Liouville theorems for fractional parabolic equations, \,\, Adv. Nolinear Stud., 21 (2021), 939-958.

\bibitem[CWu]{CWu} W. Chen, L. Wu, \,\, Uniform a priori estimates for solutions of higher critical order fractional equations,\,\, Calc. Var., 60 (2021), 102.

\bibitem[CW1]{CW1} W. Chen, L. Wu, \,\, Monotonicity and one-dimensional symmetry  of solutions for fractional reaction-diffusion equations and various applications of sliding methods,\,\, Ann. Mat. Pura Appl. 203(2024) 173-204.

\bibitem[CWW]{CWW} W. Chen, L. Wu, and P. Wang, \,\, Nonexistence of solutions for indefinite fractional  parabolic equations, \,\, Adv. Math., 392 (2021), 108018.


\bibitem[CZ]{CZ} W. Chen, J. Zhu, \,\, Indefinite fractional elliptic problem and Liouville theorems, \,\, J. Diff. Equa., 260 (2016), 4758-4785.

\bibitem[DL]{DL} Y. Du and S. Li, Nonlinear Liouville theorems and a priori
estimates for indefinite superlinear elliptic equations, Adv. Differential Equations 10 (2005), no. 8, 841-860.

\bibitem[CWX]{CWX} C. Li, Z. Wu, and H. Xu, \,\, Maximum principles and B\^{o}cher type theorems, \,\, Proc. Natl. Acad. Sci. USA 115(2018) 6976-6979.



\bibitem[DLL]{DLL}  W. Dai, Z. Liu, and G. Lu, \,\, Liouville type theorems for PDE and IE systems involving fractional Laplacian on a half space,\,\, Potential Anal., 46 (2017), 569-588.

 \bibitem[DQ]{DQ} W. Dai, G. Qin,\,\, Classification of nonnegative classical solutions to third-order equations,\,\, Adv. Math., 328 (2018), 822-857.

 \bibitem[DQ1]{DQ1} W. Dai, G. Qin, \,\,  Liouville type theorems for fractional and higher order Henon–Hardy type equations via the method of scaling spheres, Int. Math. Res. Not., 11 (2023), 9001-9070.


\bibitem[Fa]{Fa} M. Fall, \,\, Entire s-harmonic functions are affine, \,\, Proc. AMS, 144(2016) 2587-2592.

\bibitem[FC]{FC} Y. Fang and W. Chen, \,\, A Liouville type theorem for poly-harmonic Dirichlet problem in a half space, \,\, Advances in Math. 229(2012) 2835-2867.



\bibitem[Far]{Far} A. Farina, Symmetry for solutions of semilinear elliptic equations in $\mathbb{R}^N$ and related conjectures, Ricerche Math., 48 (1999), 129-154.

\bibitem[FR]{FR} X. Fern\'{a}ndez-Real, X. Ros-Oton, \,\, Regularity theory for general stable operators: parabolic equations, \,\, J. Funct. Anal.,  272 (2017), 4165-4221.

\bibitem[GLO]{GLO}, Y. Guo, C. Li, and Y. Ouyang, \,\,  Boundary regularity  and a priori estimates  for  fractional equations on unbounded domains, \,\, preprint, 2025.



\bibitem[GT]{GT}  D. Gilbarg and N. S. Trudinger,\,\,  Elliptic partial differential equations of second order. \,\,Reprint of the 1998 edition. Classics in Mathematics. Springer, Berlin, 2001.

\bibitem[GS]{GS} B. Gidas and J. Spruck, A priori bounds for positive solutions of nonlinear elliptic equations,
Comm. Partial Diff. Eqns. 8 (1981), 883-901.

\bibitem[Ha]{Ha} F. Hang, \, On the integral systems related to Hardy-Littlewood-Sobolev inequality, \, Math. Res. Lett. 14(2007) 373-383.

\bibitem[HLZ]{HLZ} X. Han, G. Lu, and J. Zhu, \,\, Characterization of balls in terms of Bessel-potential integral equation, \,\, J. Diff. Equa. 252(2012) 1589-1602.


\bibitem[JLX]{JLX} T. Jin, Y. Y. Li, and J. Xiong, \,\, On a fractional Nirenberg problem, part I: blow up analysis and compactness of solutions, \,\, J. Eur. Math. Soc. 16(2014) 1111-1171.





\bibitem[Ka]{Ka} J. Kazdan, \,\, Prescribing the curvature of a Riemannian manifold, \,\, CBMS Lecture, published by AMS, No. 57, 1984.

\bibitem[L]{L} N. S. Landkof, \,\, Foundations of modern potential theory, \,\, Springer-Verlag Berlin Heidelberg, New York, 1972. Translated from the Russian by A. P. Doohovskoy, Die Grundlehren der mathematischen Wissenschaften, Band 180.

\bibitem[Lei]{Lei} Y. Lei, \,\, Asymptotic properties of positive solutions of the Hardy-Sobolev type equations, \,\, J. Diff. Equa. 254(2013) 1774-1799.

\bibitem[LW]{LW} C. Li and L. Wu, \,\,  Pointwise regularity for fractional equations,  \,\, J. Diff. Equa., 302 (2021) 1-36.

\bibitem[LWX]{LWX} C. Li, Z. Wu, and H. Xu, Maximum principles and B\^{o}cher type theorems, \,\, Proc. Natl. Acad. Sci. USA, 115(2018) 6976-6979.

\bibitem[LZ]{LZ} G. Lu and J. Zhu, \,\, An overdetermined problem in Riesz-potential and fractional Laplacian, \,\, Nonlinear Analysis 75(2012) 3036-3048.

\bibitem[LZh]{LZh} Li, Y., Zhu, M.:{ Uniqueness theorems through the method of moving spheres}, Duke. Math. J. 80(1995) 383-417.

\bibitem[Lin]{Lin} C.S. Lin, \,\, On Liouville theorem and apriori estimates for the scalar curvature equations, \,\, Ann. Scula Norm. Sup. Pisa CI.Sci., 27(1998) 107-130.

\bibitem[LiY]{LiY} Y. Li, \,\, A semiliear Dirichlet problem involving the fractional Laplacian in $\mathbb{R}^n_+$, \,\, submitted to Adv. Nonlinear Stud., 2022.

\bibitem[LC]{LC} C. Li, W. Chen, \,\, A Hopf type lemma for fractional equations, \,\, Proc. Amer. Math. Soc., 147 (2019), 1565-1575.

\bibitem[LXYZ]{LXYZ} C. Li, M, Xu, H. Yang, and R. Zhuo, \,\, The direct moving sphere for fractional Laplace equation, \,\, J. Funct. Anal. 289(2025), 29pp.

\bibitem[MC]{MC} L. Ma, D. Chen, \,\, A Liouville type theorem for an integral system, \,\, Comm. Pure Appl. Anal., 5 (2006), 855-859.

\bibitem[MCL]{MCL}  C. Ma, W. Chen, and C. Li, \,\, Regularity of solutions for an integral system of Wolff type, \,\, Adv. Math., 3 (2011), 2676-2699.

\bibitem[MZhao]{MZhao} L. Ma, L. Zhao, \,\, Classification of positive solitary solutions of the nonlinear Choquard equation, \,\,  Arch. Ration. Mech. Anal., 195 (2010), 455-467.

\bibitem[MZ1]{MZ1} L. Ma,  Z. Zhang, \,\, Symmetry of positive solutions for Choquard equations with fractional $p$-Laplacian,\,\, Nonlinear Anal., 182 (2019), 248-262.

\bibitem[MZ2]{MZ2} L. Ma, Z. Zhang,\,\, Monotonicity for fractional Laplacian systems in unbounded Lipschitz domains, \,\,Discrete Contin. Dyn. Syst., 41 (2021), 537-552.

\bibitem[MZ3]{MZ3} L. Ma, Z. Zhang,\,\, Monotonicity of positive solutions for fractional $p$-systems in unbounded Lipschitz domains,\,\, Nonlinear Anal., 198 (2020), 111892.


\bibitem[Ni]{Ni} W.  Ni,\,\, The mathematics of diffusion,\,\,  Society for Industrial and Applied Mathematics, 2011.

\bibitem[PQ]{PQ}
P. Pol\'{a}\v{c}ik, P. Quittner, Liouville type theorems and complete blow-up for indefinite superlinear parabolic equations. Nonlinear elliptic and parabolic problems, 391-402, Progr. Nonlinear Differential Equations Appl., 64, Birkh\"{a}user, Basel, 2005.

\bibitem[PQRV]{PQRV} A. de Pablo, F. Quir\'{o}s, A. Rodr\'{i}guez, and J. L. V\'{a}zquez, \,\, A fractional porous medium equation, \,\, Adv. Math., 226 (2011), 1378-1409.

\bibitem[PQS]{PQS}
P. Pol\'{a}\v{c}ik, P. Quittner, Ph. Souplet, Singularity and decay estimates in superlinear problems via
Liouville-type theorems, Part II: parabolic equations, Indiana Univ. Math. J., 56 (2007), 879-908.

\bibitem[Q]{Q} P. Quittner, Optimal Liouville theorems for superlinear parabolic problems, arXiv:2003.13223.

\bibitem[QS]{QS}  P. Quittner and F. Simondon, A priori bounds and complete blow-up of positive solutions of indefinite superlinear parabolic problems, J. Math. Anal. Appl. 304 (2005), 614-631.





\bibitem[RS]{RS} X. Ros-Oton and J. Serra, \,\, The Dirichlet problem for the fractional Laplacian: regularity up to the boundary, \,\, J. de Math. Pures et Appl., 101(2014) 275-302.

\bibitem[RS1]{RS1} X. Ros-Oton, J. Serra, \,\, The boundary harnack principle for nonlocal elliptic operators in non-divergence form,\,\, Potential Anal., 51 (2019), 315-331.
    \bibitem[S]{S}E. M. Stein, Singular Integrals and Differentiability Properties of Functions, Princeton University Press, 1971.
\bibitem[S]{S}P.R.  Stinga,  Fractional Powers of Second Order Partial Differential Operators: Extension Problem and Regularity Theory. PhD Thesis, Universidad Autónoma de Madrid (2010).
\bibitem[Si]{Si} L. Silvestre, \,\, Regularity of the obstacle problem for a fractional power of the Laplace operator, \,\, Comm. Pure Appl. Math. 60(2007) 67-112.
\bibitem[ST]{ST} P. R. Stinga, J. L. Torrea, Regularity theory and extension problem for fractional nonlocal parabolic equations and the master equation, SIAM J. Math. Anal., 49 (2017), 3893-3924.
\bibitem[ST1]{ST1} P. R. Stinga, J. L. Torrea, Extension problem and Harnack's inequality for some fractional operators, Commun. Partial Differ. Equ., 35 (11), (2010) 2092-2122.

\bibitem[WC]{WC} P. Wang, W. Chen, \,\, Hopf's lemmas for parabolic fractional
$p$-Laplacians,\,\, Commun. Pure Appl. Anal., 21 (2022), 3055-3069.


\bibitem[WC2]{WC2} L. Wu and  W. Chen, \, The sliding methods for the fractional $p$-Laplacian, Adv. Math., 361 (2020), 106933.35.

\bibitem[WN]{WN} L. Wu and P. Niu, \, Symmetry and nonexistence of positive solutions to fractional p-Laplacian equations, \, Disc. Cont. Dyn. Sys. 39(2019), 1573-1583.


\bibitem[WuC]{WuC} L. Wu, W. Chen, \,\, Monotonicity of solutions for fractional equations with De Giorgi type nonlinearities (in Chinese),\,\,  Sci. Sin. Math., 50 (2020).

\bibitem[WuC1]{WuC1} L. Wu, W. Chen \,\, Ancient solutions to nonlocal  parabolic equations, \,\, Adv. Math, 408 (2022), 108607.

\bibitem[XFR]{XFR} X. Fern\'{a}ndez-Real, X. Ros-Oton, Regularity theory for general stable operators: parabolic equations, J. Funct. Anal.  272 (2017) 4165-4221.


\bibitem[XY]{XY} X. Xu and P. Yang, \,\,  Remarks on prescribing Gaussian curvature, \,\, Trans. AMS, 336(1993) 831-840.

\bibitem[Z]{Z} M. Zhu, Liouville theorems on some indefinite equations, Proc. Roy.
Soc. Edinburgh, 129A (1999), 649-661.

\bibitem[Z1]{Z1} M. Zhu, On elliptic problems with indefinite superlinear
boundary conditions, J. Diff. Eqns. 193 (2003), no. 1, 180-195.


\bibitem[ZCCY]{ZCCY} R. Zhuo, W. Chen, X. Cui, and Z. Yuan, \,\, Symmetry and non-existence of solutions for a nonlinear system involving the fractional Laplacian, \,\, Disc. Cont. Dyn. Sys. 36 (2016), 1125-1141.


\bibitem[ZL]{ZL} R. Zhuo, C. Li,\,\, Classification of anti-symmetric solutions to nonlinear fractional Laplace equations, \,\, Calc. Var., 61 (2022), 17.


\end{thebibliography}
\end{document}